%% file: main.tex
\documentclass{article}
\usepackage[a4paper, top=1.25in, bottom=1.25in, left=1.25in, right=1.25in]{geometry}
\usepackage[utf8]{inputenc}
\usepackage{graphicx} 

\usepackage{amsmath,amsfonts}
\usepackage{array}
\usepackage[font=footnotesize,labelfont=bf]{caption}
\usepackage{subcaption}
\usepackage{textcomp}
\usepackage{stfloats}
\usepackage{url}
\usepackage{verbatim}
\usepackage{graphicx}
\usepackage[noadjust]{cite}

\usepackage{amsmath,amsfonts,amssymb,amsthm}
\usepackage{enumitem}
\usepackage{lscape}
\usepackage{mathtools}
\usepackage{braket}
\usepackage{algorithm,algpseudocode}
\usepackage{siunitx}[=v2]
\usepackage{bm}
\usepackage{bbm}
\usepackage{mleftright}
\usepackage{booktabs}
\usepackage[usenames,dvipsnames]{xcolor}
\usepackage{pgfplots}
\usepackage{tikz}
\usepackage{color}
\usepackage{hyperref}
\hypersetup{colorlinks=true,
    linkcolor=BrickRed,
    citecolor=OliveGreen,
    urlcolor=MidnightBlue}
\usepackage{enumitem}
\usepackage{cleveref}
\usepackage{makecell}

\usepackage{afterpage}
\usepackage{booktabs}
\usepackage{graphicx}

\setlistdepth{5}
\setlist[itemize,1]{label=$\bullet$}
\setlist[itemize,2]{label=$\bullet$}
\setlist[itemize,3]{label=$\bullet$}
\setlist[itemize,4]{label=$\bullet$}
\setlist[itemize,5]{label=$\bullet$}

\renewlist{itemize}{itemize}{5}

\usepackage{listings}
\definecolor{codegreen}{rgb}{0,0.6,0}
\definecolor{codegray}{rgb}{0.5,0.5,0.5}
\definecolor{outputgray}{rgb}{0.35,0.35,0.35}
\definecolor{codepurple}{rgb}{0.58,0,0.82}
\definecolor{backcolour}{rgb}{0.95,0.95,0.92}
\lstdefinestyle{mystyle}{
    backgroundcolor=\color{backcolour},   
    commentstyle=\color{codegreen},
    keywordstyle=\color{magenta},
    numberstyle=\tiny\color{codegray},
    stringstyle=\color{codepurple},
    basicstyle=\ttfamily\footnotesize,
    breakatwhitespace=false,         
    breaklines=true,                 
    captionpos=t,                    
    keepspaces=true,                 
    numbers=left,                    
    numbersep=5pt,                  
    showspaces=false,                
    showstringspaces=false,
    showtabs=false,                  
    tabsize=2,
    frame=lines,xleftmargin=2em,framexleftmargin=1.5em
}
\pgfplotsset{compat=1.9}
\usepgfplotslibrary{fillbetween}

\DeclareMathOperator*{\minimize}{min}
\DeclareMathOperator*{\maximize}{max}
\DeclareMathOperator*{\subjto}{subj. to}

\DeclareMathOperator{\tr}{tr}

\DeclareMathOperator{\closure}{cl}
\DeclareMathOperator{\interior}{int}

\DeclareMathOperator{\domain}{dom}
\DeclareMathOperator{\diag}{diag}
\DeclareMathOperator{\image}{im}
\DeclareMathOperator{\vect}{vec}

\DeclareMathOperator{\inv}{inv}

\DeclarePairedDelimiterX{\divx}[2]{(}{)}{#1\mspace{1.5mu}\delimsize\|\mspace{1.5mu}#2}
\DeclarePairedDelimiterX{\divy}[2]{(}{)}{#1\mspace{1mu}\delimsize|\mspace{1mu}#2}
\DeclarePairedDelimiterX{\inp}[2]{\langle}{\rangle}{#1, #2}
\DeclarePairedDelimiterX{\norm}[1]{\lVert}{\rVert}{#1}
\DeclarePairedDelimiterX{\abs}[1]{\lvert}{\rvert}{#1}
\DeclarePairedDelimiterX{\bk}[2]{\langle}{\rangle}{#1 \delimsize\vert #2}

\newcommand*{\vertbar}{\rule[-1ex]{0.5pt}{2.5ex}}
\newcommand*{\horzbar}{\rule[.5ex]{2.5ex}{0.5pt}}

\makeatletter
\def\ALG@special@indent{%
    \ifdim\ALG@thistlm=0pt\relax
        \hskip-\leftmargin
    \else
        \hskip\ALG@thistlm
    \fi
}
\newcommand{\Input}[1]{\item[]\noindent\ALG@special@indent \textbf{Input:}\ #1}
\newcommand{\Output}[1]{\item[]\noindent\ALG@special@indent \textbf{Output:}\ #1}
\newcommand{\Indent}[1]{\item[]\noindent\ALG@special@indent \hspace{2.675em} #1}
\makeatother

\newcommand*\conj[1]{\overline{#1}}

\NewDocumentCommand{\grad}{e{_^}}{%
  \mathop{}\!
  \nabla
  \IfValueT{#1}{_{\!#1}}
  \IfValueT{#2}{^{#2}}
}

\makeatletter
\def\smallunderbrace#1{\mathop{\vtop{\m@th\ialign{##\crcr
   $\hfil\displaystyle{#1}\hfil$\crcr
   \noalign{\kern3\p@\nointerlineskip}%
   \tiny\upbracefill\crcr\noalign{\kern3\p@}}}}\limits}
\makeatother

\DeclareMathVersion{nxbold}
\SetSymbolFont{operators}{nxbold}{OT1}{cmr} {b}{n}
\SetSymbolFont{letters}  {nxbold}{OML}{cmm} {b}{it}
\SetSymbolFont{symbols}  {nxbold}{OMS}{cmsy}{b}{n}
\usepackage{dcolumn}
\newcolumntype{d}{D{.}{.}{-1}}
\newcommand\mc[1]{\multicolumn{1}{c}{#1}} 
\makeatletter
\newcolumntype{Z}[3]{>{\mathversion{nxbold}\DC@{#1}{#2}{#3}}c<{\DC@end}}
\makeatother

\newtheorem{thm}{Theorem}[section]
\newtheorem{lem}[thm]{Lemma}

\newtheorem{rem}[thm]{Remark}

\newcommand{\footremember}[2]{%
    \footnote{#2}
    \newcounter{#1}
    \setcounter{#1}{\value{footnote}}%
}
\newcommand{\footrecall}[1]{%
    \footnotemark[\value{#1}]%
} 

\begin{document}

\title{QICS: Quantum Information Conic Solver}

\author{%
    Kerry He\footremember{monash}{Department of Electrical and Computer Systems Engineering, Monash University, Clayton VIC 3800, Australia. \url{{kerry.he1, james.saunderson}@monash.edu}} \and James Saunderson\footrecall{monash} \and Hamza Fawzi\footremember{cambridge}{Department of Applied Mathematics and Theoretical Physics, University of Cambridge, Cambridge CB3 0WA, United Kingdom. \url{h.fawzi@damtp.cam.ac.uk}}
}
\date{}

\maketitle

\begin{abstract}
    We introduce QICS (Quantum Information Conic Solver), an open-source primal-dual interior point solver fully implemented in Python, which is focused on solving optimization problems arising in quantum information theory. QICS has the ability to solve optimization problems involving the quantum relative entropy, noncommutative perspectives of operator convex functions, and related functions. It also includes an efficient semidefinite programming solver which exploits sparsity, as well as support for Hermitian matrices. QICS is also currently supported by the Python optimization modelling software PICOS. This paper aims to document the implementation details of the algorithm and cone oracles used in QICS, and serve as a reference guide for the software. Additionally, we showcase extensive numerical experiments which demonstrate that QICS outperforms state-of-the-art quantum relative entropy programming solvers, and has comparable performance to state-of-the-art semidefinite programming solvers. 
\end{abstract}

\section{Introduction}

QICS (Quantum Information Conic Solver) is a unified conic solver which implements state-of-the-art techniques for both symmetric and nonsymmetric cone programming, and is focused on problems arising in quantum information theory. In particular, we solve conic programs of the form
\begin{gather}\label{eqn:primal}\tag{P}
    \begin{aligned}
        \minimize_{x\in\mathbb{R}^n} \quad & c^\top x \\
        \subjto \quad & b - Ax = 0 \\
                \quad & h - Gx \in \mathcal{K},
    \end{aligned}
\end{gather}
where $c\in\mathbb{R}^n$, $A\in\mathbb{R}^{p\times n}$, $b\in\mathbb{R}^p$, $G\in\mathbb{R}^{q\times n}$, $h\in\mathbb{R}^q$, and $\mathcal{K}=\mathcal{K}_1\times\ldots\times\mathcal{K}_k\subset\mathbb{R}^q$ is a Cartesian product of $k$ proper cones (i.e., convex, closed, has nonempty interior, and pointed). One can show that the corresponding conic dual problem is
\begin{gather}\label{eqn:dual}\tag{D}
    \begin{aligned}
        \maximize_{y\in\mathbb{R}^p, z\in\mathbb{R}^q} \quad & -b^\top y - h^\top z \\
        \subjto \quad & c + A^\top y + G^\top z = 0 \\
                \quad & z \in \mathcal{K}_*,
    \end{aligned}
\end{gather}
where $\mathcal{K}_*\subset\mathbb{R}^q$ is the dual cone of $\mathcal{K}$. A full list of cones currently supported by QICS can be found in Section~\ref{subsec:native-cones}. We outline our main contributions and the main features of QICS, below.

\paragraph{Efficient quantum relative entropy programming}
An important function that commonly arises in quantum information theory is the quantum (Umegaki) relative entropy, defined as
\begin{equation*}
    S\divx{X}{Y} \coloneqq \tr[X (\log(X) - \log(Y))],
\end{equation*}
for all $X\in\mathbb{H}^n_+$ and $Y\in\mathbb{H}^n_+$ satisfying $\ker(Y)\subseteq\ker(X)$ (here, we use $\mathbb{H}^n_+$ to denote the set of $n\times n$ positive semidefinite matrices, and $\ker$ to denote the kernel of a matrix). Many problems in quantum information theory involve minimizing quantum relative entropy subject to affine constraints, which we refer to as quantum relative entropy programs. These can be formulated as nonsymmetric cone programs of the form~\eqref{eqn:primal} by minimizing over the epigraph of the quantum relative entropy, i.e.,
\begin{equation*}
    \mathcal{QRE}_n \coloneqq \{ (t, X, Y) \in \mathbb{R}\times\mathbb{H}^n_{+}\times\mathbb{H}^n_{+} : t \geq S\divx{X}{Y}, \ \ker(Y)\subseteq\ker(X) \},
\end{equation*}
which we call the quantum relative entropy cone. Recently, self-concordance of the natural barrier function of the quantum relative entropy cone and other related cones was proven in~\cite{fawzi2023optimal}, which has opened up the opportunity for quantum relative entropy programs to be solved using interior point methods for nonsymmetric cone programming. Using these barriers, QICS not only provides support for optimizing over the quantum relative entropy cone, but also for the epigraph of noncommutative perspectives of operator convex functions, which include the operator relative entropy~\cite{fujii1989relative} and weighted matrix geometric means~\cite{kubo1980means}.

To solve these nonsymmetric cone programs, we use the Skajaa-Ye algorithm~\cite{nesterov2012towards} with the stepping algorithm used in Hypatia~\cite{coey2023performance}. The main advantage of the Skajaa-Ye algorithm over other interior point methods for nonsymmetric cones, e.g.,~\cite{tunccel2001generalization,myklebust2014interior,myklebust2015primal,nemirovski2005cone,karimi2020primal}, is that it does not require efficient oracles for the conjugate of the barrier function of the cone, which are currently unavailable for the quantum relative entropy cone. 

In practice, one of the main bottlenecks when solving quantum relative entropy programs with interior point methods is solving linear systems with the Hessian of the barrier function. We present two ways we can make this step easier to perform. First, we show how we can avoid this step altogether by solving an equivalent conic program involving a suitable pre-image of the cone. Second, we implement efficient oracles for some important slices of the quantum relative entropy cone introduced in~\cite{he2024exploiting}, including the quantum conditional entropy cone, as well as slices of the quantum relative entropy cone that arise when computing key rates of quantum cryptographic protocols. As shown in~\cite{he2024exploiting}, exploiting structures which arise in these slices of the quantum relative entropy cone can make solving linear systems with the Hessian of the barrier function significantly easier. Numerical experiments show that QICS is significantly faster than Hypatia~\cite{coey2023performance}, DDS~\cite{karimi2023efficient,karimi2024domain}, and CVXQUAD~\cite{fawzi2019semidefinite}, which are currently the only other available solvers and methods which can solve general quantum relative entropy programs.

\paragraph{Efficient semidefinite programming}

Interior point methods for symmetric cones have been more thoroughly studied than for general nonsymmetric cones, and robust implementations of these ideas are available in, e.g., MOSEK~\cite{mosek}, CVXOPT~\cite{vandenberghe2010cvxopt}, and SeDuMi~\cite{sturm1999using}. QICS takes advantage of some of these techniques to solve semidefinite programs. In particular, when solving symmetric cone programs, including linear, semidefinite, and second-order cone programs, we use the Nesterov-Todd algorithm~\cite{nesterov1997self,nesterov1998primal} by CVXOPT~\cite{vandenberghe2010cvxopt} instead of using the Skajaa-Ye algorithm. This algorithm has the advantage of treating the primal and dual in a symmetric way. Additionally, in practice, the Nesterov-Todd algorithm can take much more aggressive steps than the Skajaa-Ye algorithm, and therefore typically converges in significantly fewer iterations. Moreover, we exploit sparsity in the data matrices $A$ or $G$ in~\eqref{eqn:primal} when solving semidefinite programs by using the technique described in~\cite{fujisawa1997exploiting} to construct the Schur complement matrix more efficiently. Overall, numerical experiments show that our semidefinite programming solver is significantly faster than general purpose conic solvers such as CXVOPT and Hypatia, and is competitive with state-of-the-art semidefinite programming solvers including MOSEK~\cite{mosek}, SDPA~\cite{yamashita2010high,yamashita2012latest}, SDPT3~\cite{toh1999sdpt3}, and SeDuMi~\cite{sturm1999using}.

\paragraph{Complex matrices}

Often in quantum information theory, we are interested in optimization problems with complex Hermitian decision variables. Most semidefinite programming software only supports real symmetric matrices, and require users to lift Hermitian positive semidefinite cones to symmetric positive semidefinite cones of twice the size~\cite{goemans2001approximation}. A similar technique was recently shown to work for the quantum relative entropy cone~\cite{karimi2023efficient}. Instead of performing these lifts, QICS gives users the option of directly defining cones either over real symmetric or complex Hermitian matrix variables. This allows us to model the conic problem with a smaller barrier parameter, and typically results in improved computational performance in practice.

\paragraph{Benchmark libraries}

To evaluate the performance of QICS against other state-of-the-art solvers, we developed three benchmark libraries which draw from problems arising from quantum information theory. Notably, many of these problems require optimizing over Hermitian matrices. The first library of problems consists of semidefinite programs, which we store using the SDPA sparse format~\cite{fujisawa2002sdpa}. For Hermitian semidefinite programs, we slightly modify the SDPA sparse format by allowing values to be complex. The second library of problems consists of quantum relative entropy programs, and the third library of problems consists of conic programs involving the epigraphs of noncommutative perspective functions. We store these problems in the Conic Benchmark Format~\cite{friberg2016cblib}, which we have slightly modified by adding in support for all cones currently supported by QICS (both their symmetric and Hermitian counterparts, where applicable).

\subsection{Outline}
In Section~\ref{sec:software}, we introduce the QICS software, and describe how users can use QICS to solve conic programs. In Sections~\ref{sec:pdipm} and~\ref{sec:cones}, we describe some important implementation details of QICS, including the primal-dual interior point algorithm used, details about the cone oracles for the family of quantum entropy and noncommutative perspective cones we currently support, and how sparsity in the problem structure is exploited. Finally, in Section~\ref{sec:experiments}, we present numerical experiments benchmarking QICS against other state-of-the-art software for solving semidefinite, quantum relative entropy, and noncommutative perspective programs.

\subsection{Notation}

We use $\mathbb{R}^n$, $\mathbb{S}^n$, and $\mathbb{H}^n$ to denote the set of $n$-dimensional real vectors, $n\times n$ real symmetric matrices, and $n\times n$ complex Hermitian matrices, respectively. We also use $\mathbb{R}^n_+$, $\mathbb{S}^n_+$, and $\mathbb{H}^n_+$ to denote the nonnegative orthant, real positive semidefinite matrices, and complex semidefinite matrices, respectively, and $\mathbb{R}^n_{++}$, $\mathbb{S}^n_{++}$, and $\mathbb{H}^n_{++}$ to denote each of their interiors. Additionally, for Hermitian matrices $X,Y\in\mathbb{H}^n$, we use $X\succeq Y$ to mean $X-Y\in\mathbb{H}^n_{+}$, and similarly for symmetric matrices. We use $\closure$ and $\interior$ to denote the closures and interiors of a set. For a matrix $A\in\mathbb{C}^{m\times n}$, we use $A^\top$ to denotes its transpose, $\conj{A}$ to denote its conjugate, and $A^\dag$ to denote its conjugate transpose. We use $\mathbb{I}$ to denote the identity matrix, and $\{ e_i \}$ to denote the unit vectors in the standard basis.

Consider a twice differentiable function $f:\domain f\subseteq\mathbb{V}\rightarrow\mathbb{W}$, where $\mathbb{V}$ and $\mathbb{W}$ are $n$- and $m$-dimensional real vector spaces, respectively. At a point $x\in\interior\domain f$, we use $\mathsf{D}f(x)[v]\in\mathbb{W}$ to represent the directional derivative of $f$ in the direction $v\in\mathbb{V}$, and $\grad f(x) \in \mathbb{R}^{n\times m}$ to denote the gradient vector or (transposed) Jacobian of $f$. When $m=1$, we use $\mathsf{D}f(x)\in\mathbb{V}^*$ to represent the gradient of $f$ such that $\mathsf{D}f(x)[v]=\inp{v}{\mathsf{D}f(x)}$. Similarly, we use $\mathsf{D}^2f(x)[v, w]\in\mathbb{W}$ to represent the second directional derivative of $f$ in the directions $v,w\in\mathbb{V}$. When $m=1$, we use $\mathsf{D}^2f(x)[v]\in\mathbb{V}^*$ to represent the directional derivative of $\mathsf{D}f(x)$ in the direction $v$ such that $\mathsf{D}^2f(x)[v,w]=\inp{w}{\mathsf{D}^2f(x)[v]}$, and $\grad^2 f(x) \in \mathbb{R}^{n\times n}$ to denote the Hessian matrix of $f$ at $x\in\mathbb{V}$. When a function has multiple input parameters, we use subscripts to denote which variable we are taking the derivatives with respect to.

In addition to the quantum relative entropy, we define quantum (von Neumann) entropy as 
\begin{equation*}
    S(X)\coloneqq-\tr[X\log(X)],
\end{equation*}
for all $X\in\mathbb{H}^n_+$. We also define classical (Shannon) entropy as 
\begin{equation*}
    H(x)\coloneqq-\sum_{i=1}^n x_i\log(x_i),
\end{equation*}
for all $x\in\mathbb{H}^n_+$, and classical relative entropy (Kullback-Leibler divergence) as 
\begin{equation*}
    H\divx{x}{y}\coloneqq\sum_{i=1}^nx_i\log(x_i/y_i),
\end{equation*}
for all $x\in\mathbb{R}^n_+$ and $y\in\mathbb{R}^n_+$ satisfying $\ker(\diag(y))\subseteq\ker(\diag(x))$. Let $I\subseteq\mathbb{R}$ be an interval. For a scalar function $g:I\subseteq\mathbb{R}\rightarrow\mathbb{R}$, for each positive integer $n$ we define its noncommutative perspective $P_g:\domain P_g\subseteq \mathbb{H}^n\times\mathbb{H}^n\rightarrow\mathbb{H}$ by 
\begin{equation*}
    P_g(X, Y) \coloneqq X^{1/2} g(X^{-1/2} Y X^{-1/2}) X^{1/2},
\end{equation*}
where $\domain P_g = \{X\in\mathbb{H}^n: \lambda(X)\in I\}$ and $\lambda(X)$ denotes the eigenvalues of $X$. Consider a system composed of $r$ subsystems, where the $i$-th subsystem has dimension $n_i$, and let $n=\Pi_i n_i$. The partial trace $\tr_i:\mathbb{H}^n\rightarrow\mathbb{H}^{n/n_i}$ is defined as the unique linear operator satisfying
\begin{equation*}
    \tr_i(X_1 \otimes \ldots \otimes X_{i} \otimes \ldots \otimes X_{r}) = X_1 \otimes \ldots \otimes \tr[X_{i}] \otimes \ldots \otimes X_{r},
\end{equation*}
for all $X_j\in\mathbb{H}^{n_j}$ and $j=1,\ldots,r$, where $\otimes$ denotes the Kronecker product. Finally, we denote $\odot$ as the elementwise (Hadamard) multiplication, and $\oslash$ as the elementwise division between two matrices.

\section{Software interface}\label{sec:software}

The algorithms and cone oracles discussed in Sections~\ref{sec:pdipm} and~\ref{sec:cones} are implemented in the Python package QICS, which is currently supported for Python versions 3.8 or later. The only dependencies for QICS are NumPy, SciPy, and Numba. QICS can be directly installed from the Python Package Index by using
\begin{center}
    \texttt{\$ pip install qics}
\end{center}
Alternatively, our source code can be found at
\begin{center}
    \url{https://github.com/kerry-he/qics}.
\end{center}
In Section~\ref{subsec:native}, we describe the native interface QICS provides for users to define conic problems. In Section~\ref{subsec:picos}, we explain how users can use the high-level optimization modelling interface PICOS~\cite{picos} to more easily parse problems and solve them using QICS.

\subsection{Native interface}\label{subsec:native}

There are three main classes users interact with when defining a conic problem using the native interface, which we introduce in \Cref{subsec:native-cones,subsec:native-model,subsec:native-solver}. In Section~\ref{subsec:native-matrix}, we discuss how symmetric and Hermitian matrix variables are represented in QICS. Finally, in Section~\ref{subsec:native-example}, we walk through an example for how to express a simple quantum relative entropy program into a form that is appropriate for QICS.

\subsubsection{Cones}\label{subsec:native-cones}
Users define the Cartesian product of cones $\mathcal{K}$ by defining a list of cone classes from the \texttt{qics.cones} module. Below, we list the definitions and interfaces to all of the cones QICS currently supports. Implementation details for the family of quantum entropy and noncommutative perspective cones can be found in Section~\ref{sec:cones}. Note that the \texttt{iscomplex} parameter is always \texttt{False} if unspecified.

\newpage
\begin{itemize}[leftmargin=0em]
    \item[] \textbf{Symmetric cones}
    \begin{itemize}
        \item[$\bullet$] \texttt{qics.cones.NonNegOrthant(n)} represents the nonnegative orthant
        \begin{equation*}
            \mathbb{R}^n_+ = \{ x \in \mathbb{R}^n : x \geq 0 \},
        \end{equation*}
        where \texttt{n} is a positive integer defining the dimension of the cone.
        \item[$\bullet$] \texttt{qics.cones.PosSemidefinite(n, iscomplex)} represents the positive semidefinite cone
        \begin{equation*}
            \mathbb{H}^n_+ = \{ X \in \mathbb{H}^n : X \succeq 0 \},
        \end{equation*}
        where \texttt{n} is a positive integer, and $X$ is defined over $\mathbb{H}^n$ when \texttt{iscomplex=True}, or restricted to $\mathbb{S}^n$ when \texttt{iscomplex=False}.
        \item[$\bullet$] \texttt{qics.cones.SecondOrder(n)} represents the second order cone
        \begin{equation*}
            \mathcal{Q}_{n+1} = \{ (t, x) \in \mathbb{R} \times \mathbb{R}^{n} : t \geq \norm{x}_2 \},
        \end{equation*}
        where \texttt{n} is a positive integer defining the dimension of the variable $x$.
    \end{itemize}

    \item[] \textbf{Classical entropy cones}
    \begin{itemize}
        \item[$\bullet$] \texttt{qics.cones.ClassEntr(n)} represents the (homogenized) classical entropy cone
        \begin{equation*}
            \mathcal{CE}_n = \closure\{ (t, u, x) \in \mathbb{R} \times \mathbb{R}_{++} \times \mathbb{R}^n_{++} : t \geq -u H(u^{-1} x) \},
        \end{equation*}
        where \texttt{n} is a positive integer defining the dimension of the variable $x$.
        \item[$\bullet$] \texttt{qics.cones.ClassRelEntr(n)} represents the classical relative entropy cone
        \begin{equation*}
            \mathcal{CRE}_n = \closure\{ (t, x, y) \in \mathbb{R} \times \mathbb{R}^n_{++} \times \mathbb{R}^n_{++} : t \geq H\divx{x}{y} \},
        \end{equation*}
        where \texttt{n} is a positive integer defining the dimension of the variables $x$ and $y$.
    \end{itemize}

    \item[] \textbf{Quantum entropy cones}
    \begin{itemize}
        \item[$\bullet$] \texttt{qics.cones.QuantEntr(n, iscomplex)} represents the (homogenized) quantum entropy cone
        \begin{equation*}
            \mathcal{QE}_n = \closure\{ (t, u, X) \in \mathbb{R} \times \mathbb{R}_{++} \times \mathbb{H}^n_{++} : t \geq -u S(u^{-1} X) \},
        \end{equation*}
        where \texttt{n} is a positive integer, and $X$ is defined over $\mathbb{H}^n$ when \texttt{iscomplex=True}, or restricted to $\mathbb{S}^n$ when \texttt{iscomplex=False}.
        \item[$\bullet$] \texttt{qics.cones.QuantRelEntr(n, iscomplex)} represents the quantum relative entropy cone
        \begin{equation*}
            \mathcal{QRE}_n = \closure\{ (t, X, Y) \in \mathbb{R} \times \mathbb{H}^n_{++} \times \mathbb{H}^n_{++} : t \geq S\divx{X}{Y} \},
        \end{equation*}
        where \texttt{n} is a positive integer, and $X$ and $Y$ are defined over $\mathbb{H}^n$ when \texttt{iscomplex=True}, or restricted to $\mathbb{S}^n$ when \texttt{iscomplex=False}.
        \item[$\bullet$] \texttt{qics.cones.QuantCondEntr(dims, sys, iscomplex)} represents the quantum conditional entropy cone
        \begin{equation*}
            \mathcal{QCE}_{\{n_i\}, j} = \closure\{ (t, X) \in \mathbb{R} \times \mathbb{H}^{\Pi_in_i}_{++} : t \geq -S(X) + S(\text{tr}_i(X)) \},
        \end{equation*}
        where 
        \begin{itemize}
            \item \texttt{dims} is a list of positive integers representing the dimensions $\{n_i\}$ of the subsystems which $X$ is defined over,
            \item  \texttt{sys} is an integer or list of integers representing which subsystems $i$ are traced out,
            \item $X$ is defined over $\mathbb{H}^{\Pi_in_i}$ when \texttt{iscomplex=True}, or restricted to $\mathbb{S}^{\Pi_in_i}$ when \texttt{iscomplex} \texttt{=False}.
        \end{itemize}
        \item[$\bullet$] \sloppy \texttt{qics.cones.QuantKeyDist(G\_info, Z\_info, iscomplex)} represents the quantum key distribution cone
        \begin{equation*}
            \mathcal{QKD}_{\mathcal{G},\mathcal{Z}} = \closure\{ (t, X) \in \mathbb{R} \times \mathbb{H}^n_{++} : t \geq -S(\mathcal{G}(X)) + S(\mathcal{Z}(\mathcal{G}(X))) \},
        \end{equation*}
        which arises in quantum key cryptography applications when aiming to certify the security of a given protocol~\cite{coles2016numerical,winick2018reliable,hu2022robust}. 
        \begin{itemize}
            \item \texttt{G\_info} specifies the linear map $\mathcal{G}$. There are two ways to specify this linear map.
            \begin{itemize}
                \item If \texttt{G\_info} is an integer, then $\mathcal{G}(X)=X$ and \texttt{G\_info} specifies the dimension of $X$.
                \item If \texttt{G\_info} is a list of NumPy arrays, then \texttt{G\_info} specifies the Kraus operators $K_i$ such that
                \begin{equation*}
                    \mathcal{G}(X) = \sum_{i} K_i X K_i^\dagger.
                \end{equation*}
            \end{itemize}
            \item \texttt{Z\_info} specifies the pinching map $\mathcal{Z}$, which has the form
            \begin{equation*}
                \mathcal{Z}(X) = \sum_{i} Z_i X Z_i^\dagger.
            \end{equation*}
            There are three ways to specify this linear map.
            \begin{itemize}
                \item If \texttt{Z\_info=r} is a positive integer, then $Z_i=e_ie_i^\top \otimes \mathbb{I}$ for $i=1,\ldots,r$.
                \item If \texttt{Z\_info=(dims, sys)} is a tuple where \texttt{dims=(n0, n1)} is a list of two positive integers and \texttt{sys} is an integer, then 
                \begin{itemize}
                    \item $Z_i=e_ie_i^\top \otimes \mathbb{I}_{n_1}$ for $i=1,\ldots,n_0$ if \texttt{sys=0}, and
                    \item $Z_i=\mathbb{I}_{n_0}\otimes e_ie_i^\top$ for $i=1,\ldots,n_1$ if \texttt{sys=1}.
                \end{itemize}
                We generalize this definition to the case where \texttt{dims} and \texttt{sys} are lists of any length.
                \item If \texttt{Z\_info} is a list of NumPy arrays, then \texttt{Z\_info} directly specifies the Kraus oeprator $Z_i$. Note that these Kraus operators must have a similar structure to those defined using the other options, i.e., must be diagonal matrices consisting of either ones and zeros, and $Z_iZ_j=0$ for all $i\neq j$.
            \end{itemize}
            \item $X$ is defined over $\mathbb{H}^n$ when \texttt{iscomplex=True}, or restricted to $\mathbb{S}^n$ when \texttt{iscomplex=False}.
        \end{itemize}
    \end{itemize}

    \item[] \textbf{Noncommutative perspective cones}
    \begin{itemize}
        \item[$\bullet$] \texttt{qics.cones.OpPerspecEpi(n, func, iscomplex)} represents the epigraph of the noncommutative perspective
        \begin{equation*}
            \mathcal{OPE}_{n, g} = \closure\{ (T, X, Y) \in \mathbb{H}^n \times \mathbb{H}^n_{++} \times \mathbb{H}^n_{++} : T \succeq P_g(X, Y) \},
        \end{equation*}
        where \texttt{n} is a positive integer, \texttt{func} determines the function $g$ (see below), and $T$, $X$, and $Y$ are defined over $\mathbb{H}^n$ when \texttt{iscomplex=True}, or restricted to $\mathbb{S}^n$ when \texttt{iscomplex=False}.
        \item[$\bullet$] \sloppy \texttt{qics.cones.OpPerspecTr(n, func, iscomplex)} represents the trace noncommutative perspective cone
        \begin{equation*}
            \mathcal{OPT}_{n, g} = \closure\{ (t, X, Y) \in \mathbb{R} \times \mathbb{H}^n_{++} \times \mathbb{H}^n_{++} : t \geq \text{tr}[P_g(X, Y)] \},
        \end{equation*}
        where \texttt{n} is a positive integer, \texttt{func} determines the function $g$ (see below), and $X$ and $Y$ are defined over $\mathbb{H}^n$ when \texttt{iscomplex=True}, or restricted to $\mathbb{S}^n$ when \texttt{iscomplex=False}.
    \end{itemize}
    For both of these noncommutative perspective cones, \texttt{func} can be defined in the following ways:
    \begin{itemize}
        \item $g(x)=-\log(x)$ if \texttt{func="log"},
        \item $g(x)=-x^p$ if \texttt{func=p} is a float where $p\in(0, 1)$, and
        \item $g(x)=x^p$ if \texttt{func=p} is a float where $p\in[-1, 0)\cup(1, 2)$.
    \end{itemize}
    
\end{itemize}

\subsubsection{Modelling}\label{subsec:native-model}
A conic program~\eqref{eqn:primal} is represented using the \texttt{qics.Model} class, and can be initialized by calling
\begin{center}
    \texttt{model = qics.Model(c, A, b, G, h, cones, offset)}
\end{center}
where
\begin{itemize}
    \item \texttt{c} is an $n \times 1$ NumPy float array representing the linear objective $c$.
    \item \texttt{A} is a $p \times n$ NumPy or sparse SciPy float array representing the linear equality constraint matrix $A$.
    \item \texttt{b} is a $p \times 1$ NumPy float array representing the linear equality constraint vector $b$.
    \item \texttt{G} is a $q \times n$ NumPy or sparse SciPy float array representing the linear cone constraint matrix $G$.
    \item \texttt{h} is a $p \times n$ NumPy float array representing the linear cone constraint vector $h$.
    \item \texttt{cones} is a list of \texttt{qics.cones} classes representing the Cartesian product of cones $\mathcal{K}$.
    \item \texttt{offset} is an optional argument specifying a constant objective offset. Default is \texttt{0}. 
\end{itemize}
A list of all currently supported \texttt{qics.cones} can be found in Section~\ref{subsec:native-cones}. When the arguments \texttt{G} and \texttt{h} are unspecified (or set to \texttt{None}), then QICS solves the simplified conic program
\begin{equation*}
    \minimize_{x\in\mathbb{R}^n} \quad c^\top x \quad \subjto \quad b - Ax = 0, \quad x \in \mathcal{K},
\end{equation*}
which uses a simplified method to solve the Newton equations (see Section~\ref{sec:newton}).

\subsubsection{Solving}\label{subsec:native-solver}
The final class users interface with is \texttt{qics.Solver}, which is used to solve the conic program defined by a \texttt{qics.Model}. The interface for this class is
\begin{center}
    \texttt{solver = qics.Solver(model, **kwargs)}
\end{center}
where \texttt{model} is an instance of a \texttt{qics.Model}, and optional solver settings include
\begin{itemize}
    \item \texttt{max\_iter} is the maximum number of solver iterations before terminating. Default is $100$.
    \item \texttt{max\_time} is the maximum time elapsed, in seconds, before terminating. Default is $3600$.
    \item \texttt{tol\_gap} is the stopping tolerance $\varepsilon_g$ for the optimality gap. Default is $10^{-8}$.
    \item \texttt{tol\_feas} is the stopping tolerance $\varepsilon_f$ for the primal-dual feasibility. Default is $10^{-8}$.
    \item \texttt{tol\_infeas} is the tolerance $\varepsilon_i$ for detecting an infeasible problem. Default is $10^{-12}$.
    \item \texttt{tol\_ip} is the tolerance $\varepsilon_{ip}$ for detecting an ill-posed problem. Default is $10^{-13}$.
    \item \texttt{tol\_near} is the multiplicative margin for certifying near optimality or infeasiblity. Default is $1000$.
    \item \texttt{verbose} is the verbosity level of the solver, where
    \begin{itemize}
        \item \texttt{verbose=0} returns no console output,
        \item \texttt{verbose=1} only prints a summary of the problem input and solution,
        \item \texttt{verbose=2} also prints a summary of the solver status at each iteration, and
        \item \texttt{verbose=3} also prints a summary of the stepper status at each iteration.
    \end{itemize}
    Default is \texttt{2}.
    \item \texttt{ir} is a Boolean which specifies whether to use iterative refinement when solving the Newton system. Default is \texttt{True}.
    \item \texttt{toa} is a Boolean which specifies whether to use third-order adjustments to improve the stepping directions. Default is \texttt{True}.    
    \item \texttt{init\_pnt} is a \texttt{qics.Point} which specifies where to initialize the interior-point algorithm from. Default is \texttt{None}, which uses the initialization strategy described in Section~\ref{eqn:initialization}.    
    \item \texttt{use\_invhess} is a Boolean which specifies whether or not to avoid using inverse Hessian-vector products by using the technique described in Section~\ref{eqn:avoid-inv-hess}. By default, this is \texttt{False} if $G\neq-\mathbb{I}$, $G$ is full column rank, and $\mathcal{K}$ mainly consists of quantum relative entropy or noncommutative perspective cones, and is \texttt{True} otherwise.
\end{itemize}
See Section~\ref{sec:termination} for details about the termination criteria. Once a \texttt{qics.Solver} has been initialized, the conic program can be solved by calling
\begin{center}
    \texttt{info = solver.solve()}
\end{center}
This returns a dictionary \texttt{info} which summarizes the solution of the conic program, and has the following keys
\begin{itemize}
    \item \texttt{x\_opt}, \texttt{y\_opt}, \texttt{z\_opt}, and \texttt{s\_opt} are the optimal variables $x^*$, $y^*$, $z^*$, and $s^*=h-Gx^*$ of the original primal~\eqref{eqn:primal} and dual~\eqref{eqn:dual} conic programs.
    \item \texttt{sol\_status} and \texttt{exit\_status} are strings describing the solution status and exit status of the solver, respectively.
    \item \texttt{num\_iter} and \texttt{solve\_time} are the total number of solver iterations and time elapsed in seconds, respectively, before the solver terminated.
    \item \texttt{p\_obj} and \texttt{d\_obj} are the primal and dual objective values, and \texttt{opt\_gap} is the (relative) optimality gap.
    \item \texttt{p\_feas} and \texttt{d\_feas} are the (relative) residuals of the primal and dual constraints.
\end{itemize}


\subsubsection{Representing matrices}\label{subsec:native-matrix}

To supply a conic program model to QICS, users are required to express the linear constraints $A$ and $G$ in matrix form. This is straightforward when variables correspond to real vectors. For variables corresponding to symmetric or Hermitian matrices, we need to first vectorize these matrices. For a real symmetric matrix $X \in \mathbb{S}^n$, we vectorize this matrix by stacking the rows of the matrix side-by-side, i.e.,
\begin{equation*}
\text{vec}(X) = \text{vec}\left(\begin{bmatrix}
                             \horzbar & x_{1}^\top & \horzbar \\
                             \horzbar & x_{2}^\top & \horzbar \\
                                      & \vdots      &         \\
                             \horzbar & x_{n}^\top & \horzbar
                         \end{bmatrix}\right) = \begin{bmatrix}
                                                     x_{1}  \\
                                                     x_{2}  \\
                                                     \vdots \\
                                                     x_{n}
                                                 \end{bmatrix} \in \mathbb{R}^{n^2}.   
\end{equation*}
This operation can be performed by using 
\begin{center}
    \texttt{x = qics.vectorize.mat\_to\_vec(X)}
\end{center}
where \texttt{X} is a NumPy array with real-valued entries.
Now instead consider a complex Hermitian matrix $X \in \mathbb{H}^n$. To vectorize this matrix, we first convert each row $x_i\in\mathbb{C}^n$ of $X$ into a real vector by splitting the real and imaginary components of each entry as follows
\begin{equation*}
    \text{split}(x_i) = \text{split}\!\left( \begin{bmatrix}
           x_{i1} \\
           x_{i2} \\
           \vdots \\
           x_{in}
        \end{bmatrix} \right) = \begin{bmatrix}
                                   \text{Re}(x_{i1}) \\
                                   \text{Im}(x_{i1}) \\
                                   \text{Re}(x_{i2}) \\
                                   \text{Im}(x_{i2}) \\
                                   \vdots \\
                                   \text{Re}(x_{in}) \\
                                   \text{Im}(x_{in}) \\
                                \end{bmatrix} \in \mathbb{R}^{2n}.
\end{equation*}
Using this, Hermitian matrices are vectorized as follows
\begin{equation*}
    \text{vec}(X) = \text{vec}\!\left(\begin{bmatrix}
             \horzbar & x_{1}^\top & \horzbar \\
             \horzbar & x_{2}^\top & \horzbar \\
                      & \vdots      &         \\
             \horzbar & x_{3}^\top & \horzbar
         \end{bmatrix}\right) = \begin{bmatrix}
                                     \text{split}(x_{1}) \\
                                     \text{split}(x_{2}) \\
                                     \vdots \\
                                     \text{split}(x_{n})
                                 \end{bmatrix} \in\mathbb{R}^{2n^2}.
\end{equation*}
In practice, we can perform this operation by using the same \texttt{qics.vectorize.mat\_to\_vec} function as for the real symmetric case, except when the argument \texttt{X} is a NumPy array with complex-valued entries.
To see how we use can use these vectorizations to represent linear constraints, consider constraints of the form 
\begin{equation*}
    \tr[A_i X] = b_i, \qquad \forall\ i=1,\ldots,p,
\end{equation*}
where $X \in \mathbb{H}^n$ is our matrix variable, and $A_i \in \mathbb{H}^n$ and $b_i\in\mathbb{R}$ encode linear constraints for $i=1,\ldots,p$. We can represent these constraints in matrix form as
\begin{equation*}
    A\text{vec}(X) = b, \qquad \text{where} \qquad A = \begin{bmatrix}
             \horzbar & \text{vec}(A_1)^\top & \horzbar \\
             \horzbar & \text{vec}(A_2)^\top & \horzbar \\
                      & \vdots               &          \\
             \horzbar & \text{vec}(A_p)^\top & \horzbar
         \end{bmatrix} \in \mathbb{R}^{p\times n^2}.
\end{equation*}
Alternatively, consider constraints of the form
\begin{equation*}
    H - \sum_{i=1}^q x_i G_i \in \mathcal{K},
\end{equation*}
where $x \in \mathbb{R}^q$ is our variable, and $H \in \mathbb{H}^n$ and $G_i \in \mathbb{H}^n$ encode linear cone constraints for $i=1,\ldots,q$. Then we can represent this in matrix form as
\begin{equation*}
    \text{vec}(H) - Gx \in \mathcal{K}, \qquad \text{where} \qquad G = \begin{bmatrix}
             \vertbar & \vertbar &        & \vertbar \\
             \text{vec}(G_1)         & \text{vec}(G_2)         & \cdots & \text{vec}(G_q) \\
             \vertbar & \vertbar &        & \vertbar
         \end{bmatrix} \in \mathbb{R}^{n^2\times q}.
\end{equation*}

\subsubsection{Illustrative example}\label{subsec:native-example}

To illustrate how to use QICS, we show how we can solve the simple quantum relative entropy program
\begin{equation}\label{eqn:example-qrep}
    \min_{Y \in \mathbb{S}^2} \quad S( X \| Y ) \quad \subjto \quad Y_{11} = Y_{22} = 1, \ Y \succeq 0,
\end{equation}
where 
\begin{equation*}
    X = \begin{bmatrix} 2 & 1 \\ 1 & 2 \end{bmatrix}.
\end{equation*}
First, we need to reformulate the above problem as a standard form conic program~\eqref{eqn:primal}. We can do this by rewriting the problem as
\begin{align*}
    \min_{t, X, Y} \quad & t \\
    \subjto \quad & \langle A_i, X \rangle = a_i, \quad i=1,2,3\\
    & \langle B_j, Y \rangle = b_j, \quad j=1,2 \\
    & (t, X, Y) \in \mathcal{QRE}_2,  
\end{align*}
where $a_1=a_2=a_3=2$, $b_1=b_2=1$, and
\begin{equation*}
    A_1 = \begin{bmatrix} 1 & 0 \\ 0 & 0 \end{bmatrix}, \quad
    A_2 = \begin{bmatrix} 0 & 1 \\ 1 & 0 \end{bmatrix}, \quad
    A_3 = \begin{bmatrix} 0 & 0 \\ 0 & 1 \end{bmatrix}, \quad
    B_1 = \begin{bmatrix} 1 & 0 \\ 0 & 0 \end{bmatrix}, \quad
    B_2 = \begin{bmatrix} 0 & 0 \\ 0 & 1 \end{bmatrix}.
\end{equation*}
Following our discussion in Section~\ref{subsec:native-matrix}, we represent the decision variable $(t, X, Y)\in\mathbb{R}\times\mathbb{S}^2\times\mathbb{S}^2$ as a vector $x\in\mathbb{R}^9$ with elements represented by
\begin{align*}
    x &= \begin{bmatrix} t & \vect(X)^\top & \vect(Y)^\top \end{bmatrix}^\top\\
      &= \begin{bmatrix} t & X_{11} & X_{12} & X_{21} & X_{22} & Y_{11} & Y_{12} & Y_{21} & Y_{22} \end{bmatrix}^\top.
\end{align*}
Given this representation of the decision variable, we can now represent our linear objective function as $c^\top x$, where
\begin{equation*}
    c = \begin{bmatrix} 1 & 0 & 0 & 0 & 0 & 0 & 0 & 0 & 0 \end{bmatrix}^\top.
\end{equation*}
We represent this objective function using a NumPy array.
\begin{lstlisting}[language=python]
import numpy
c = numpy.array([[1., 0., 0., 0., 0., 0., 0., 0., 0.]]).T
\end{lstlisting}
Additionaly, we can represent our linear equality constraints as $Ax=b$, where
\begin{equation*}
    A = \begin{bmatrix}
        0 & \vect(A_1)^\top & \vect(0)^\top \\
        0 & \vect(A_2)^\top & \vect(0)^\top \\
        0 & \vect(A_3)^\top & \vect(0)^\top \\
        0 & \vect(0)^\top & \vect(B_1)^\top \\
        0 & \vect(0)^\top & \vect(B_2)^\top
    \end{bmatrix} = \begin{bmatrix}
        0 & 1 & 0 & 0 & 0 & 0 & 0 & 0 & 0 \\
        0 & 0 & 1 & 1 & 0 & 0 & 0 & 0 & 0 \\
        0 & 0 & 0 & 0 & 1 & 0 & 0 & 0 & 0 \\
        0 & 0 & 0 & 0 & 0 & 1 & 0 & 0 & 0 \\
        0 & 0 & 0 & 0 & 0 & 0 & 0 & 0 & 1
    \end{bmatrix},
\end{equation*}
and
\begin{equation*}
    b = \begin{bmatrix} 2 & 2 & 2 & 1 & 1 \end{bmatrix}^\top.
\end{equation*}
Again, we represent these linear constraints using NumPy arrays. 
\begin{lstlisting}[language=python, firstnumber=3]
A = numpy.array([                        
    [0., 1., 0., 0., 0., 0., 0., 0., 0.],
    [0., 0., 1., 1., 0., 0., 0., 0., 0.],
    [0., 0., 0., 0., 1., 0., 0., 0., 0.],
    [0., 0., 0., 0., 0., 1., 0., 0., 0.],
    [0., 0., 0., 0., 0., 0., 0., 0., 1.] 
])
b = numpy.array([[2., 2., 2., 1., 1.]]).T
\end{lstlisting}
Next, we need to tell QICS that $x$ must be constrained within the quantum relative entropy cone $\mathcal{QRE}_2$. We do this by initializing a \texttt{qics.cones.QuantRelEntr} object.
\begin{lstlisting}[language=python, firstnumber=11]
import qics
cones = [qics.cones.QuantRelEntr(2)]
\end{lstlisting}
Finally, we initialize a \texttt{qics.Model} object to represent our conic program, and a \texttt{qics.Solver} object to solve the conic program with default settings.
\begin{lstlisting}[language=python, firstnumber=13]
model = qics.Model(c=c, A=A, b=b, cones=cones)
solver = qics.Solver(model)
info = solver.solve()
\end{lstlisting}
We can access the optimal variable $Y$ through the solver output \texttt{info} by calling
\begin{lstlisting}[language=python, firstnumber=16]
print("Optimal matrix variable Y is:")
print(info["s_opt"][0][2])
\end{lstlisting}
Running the above code using the default setting for \texttt{verbose} yields the following output.
\begin{lstlisting}[numbers=none,basicstyle=\color{outputgray}\ttfamily\footnotesize]
====================================================================
           QICS v1.0.0 - Quantum Information Conic Solver
              by K. He, J. Saunderson, H. Fawzi (2024)
====================================================================
Problem summary:
        no. vars:     9                         barr. par:    6
        no. constr:   5                         symmetric:    False
        cone dim:     9                         complex:      False
        no. cones:    1                         sparse:       False

...

Solution summary
        sol. status:  optimal                   num. iter:    7
        exit status:  solved                    solve time:   0.534
        primal obj:   2.772588704578e+00        primal feas:  6.28e-09
        dual obj:     2.772588709102e+00        dual feas:    3.14e-09

Optimal matrix variable Y is:
[[1.  0.5]
 [0.5 1. ]]
\end{lstlisting}

\subsection{PICOS interface}\label{subsec:picos}

Alternatively, users can use QICS through PICOS~\cite{picos}, a high-level Python interface for users to parse convex optimization problems to a solver. Notably, with the release of QICS, PICOS now supports a range of quantum entropy and noncommutative perspective expressions summarized in Table~\ref{tab:picos}.
Scalar functions (i.e., quantum entropy, quantum relative entropy, quantum conditional entropy, and quantum key distribution) can be used by either incorporating them in the objective function, e.g.,
\begin{center}
    \texttt{P.set\_objective("min", picos.quantrelentr(X, Y))}
\end{center}
or as an inequality constraint, e.g.,
\begin{center}
    \texttt{P.add\_constraint(t > picos.quantrelentr(X, Y))}
\end{center}
Matrix-valued functions (i.e., the operator relative entropy and matrix geometric mean) can be used in a matrix inequality expression, e.g.,
\begin{center}
    \texttt{P.add\_constraint(T >> picos.oprelentr(X, Y))}
\end{center}
or composed with a trace function to represent the corresponding scalar valued function, e.g.,
\begin{center}
    \texttt{P.set\_objective("min", picos.trace(picos.oprelentr(X, Y)))}
\end{center}
Note that these objectives and constraints need to define a convex optimization problem. Once a PICOS problem has been defined, it can be solved using QICS by calling
\begin{center}
    \texttt{P.solve(solver="qics")}
\end{center}
Below, we illustrate how the same simple quantum relative entropy program~\eqref{eqn:example-qrep} discussed in the previous section can be solved using PICOS. 
\begin{lstlisting}[language=python]
import picos

# Define the conic program
P = picos.Problem()
X = picos.Constant("X", [[2., 1.], [1., 2.]])
Y = picos.SymmetricVariable("Y", 2)

P.set_objective("min", picos.quantrelentr(X, Y))
P.add_constraint(picos.maindiag(Y) == 1)

# Solve the conic program
P.solve(solver="qics")
print("Optimal matrix variable Y is: \n", Y)
\end{lstlisting}
\begin{lstlisting}[numbers=none,basicstyle=\color{outputgray}\ttfamily\footnotesize]
Optimal matrix variable Y is:  
[ 1.00e+00  5.00e-01]
[ 5.00e-01  1.00e+00]
\end{lstlisting}

\begin{table}[t]
\centering
\caption{List of quantum entropic and noncommutative perspective functions supported by PICOS. These can either be used in the objective function (for scalar functions), or in inequality constraints to represent a convex conic problem to be solved with QICS.}
\label{tab:picos}
\resizebox{\textwidth}{!}{%
\begin{tabular}{@{}lccc@{}}
\toprule
\textbf{Function} & \textbf{PICOS expression} & \textbf{Convexity} & \textbf{Description} \\ \midrule
Quantum entropy & \texttt{picos.quantentr} & Concave & $S(X)=-\tr[X\log(X)]$\\ \addlinespace
\makecell[l]{Quantum\\relative entropy} & \texttt{picos.quantrelentr} & Convex & $S\divx{X}{Y}=\tr[X\log(X)-X\log(Y)]$ \\ \addlinespace
\makecell[l]{Quantum\\conditional entropy} & \texttt{picos.quantcondentr} & Concave & $S(X) - S(\tr_i(X))$ \\ \addlinespace
\makecell[l]{Quantum\\key distribution} & \texttt{picos.quantkeydist} & Convex & $-S(\mathcal{G}(X)) + S(\mathcal{G}(\mathcal{Z}(X))) $ \\ \addlinespace
\makecell[l]{Operator\\relative entropy} & \texttt{picos.oprelentr} & Operator convex & $-P_{\log}(X, Y) = X^{1/2} \log(X^{1/2} Y^{-1} X^{1/2}) X^{1/2}$ \\ \addlinespace
\makecell[l]{Matrix\\geometric mean} & \texttt{picos.mtxgeomean} & \makecell{Operator convex\\if $t\in[-1, 0]\cup[1, 2]$\\Operator concave\\if $t\in[0, 1]$} & $X \#_t Y = X^{1/2} (X^{-1/2} Y^{1} X^{-1/2})^t X^{1/2}$ \\ \bottomrule
\end{tabular}
}
\end{table}

\section{Primal-dual interior point algorithm}\label{sec:pdipm}

To solve the pair of primal~\eqref{eqn:primal} and dual~\eqref{eqn:dual} conic programs, we will use the homogeneous self-dual embedding of the optimality conditions by~\cite{vandenberghe2010cvxopt}, i.e., we want to find $\omega\coloneqq(x, y, z, s, \tau, \kappa) \in \mathbb{R}^n\times\mathbb{R}^p\times\mathbb{R}^q\times\mathbb{R}^q\times\mathbb{R}\times\mathbb{R}$ which satisfies
\begin{equation}\label{eqn:hsde}
    \begin{gathered}
        L(\omega) \coloneqq \begin{bmatrix}
            0 & A^\top & G^\top & c \\
            -A & 0 & 0 & b \\
            -G & 0 & 0 & h \\
            -c^\top & -b^\top & -h^\top & 0
        \end{bmatrix} \begin{bmatrix}
            x \\ y \\ z \\ \tau
        \end{bmatrix} - \begin{bmatrix}
            0 \\ 0 \\ s \\ \kappa
        \end{bmatrix} = 0 \\
        (s, \kappa, z, \tau) \in \mathcal{K} \times \mathbb{R}_+ \times \mathcal{K}_* \times \mathbb{R}_+.
    \end{gathered}
\end{equation}
The homogeneous self-dual embedding is always feasible (e.g., when $\omega=0$), while also encoding the feasiblity of the original problem through the parameters $\tau$ and $\kappa$. In particular, for a solution $\omega^*$ of~\eqref{eqn:hsde} there are three scenarios that can occur
\begin{itemize}
    \item If $\tau^*>0$ and $\kappa^*=0$, then $x^*/\tau^*$ is a solution to~\eqref{eqn:primal}, and $(y^*, z^*)/\tau^*$ is a solution to~\eqref{eqn:dual}.
    \item If $\tau^*=0$ and $\kappa^*>0$, then if $b^\top y + h^\top z < 0$ then~\eqref{eqn:primal} is infeasible, and if $c^\top x < 0$ then~\eqref{eqn:dual} is infeasible.
    \item If $\tau^*=0$ and $\kappa^*=0$, then no conclusion can be made about the optimal values or feasibility of~\eqref{eqn:primal} and~\eqref{eqn:dual}.
\end{itemize}
Recalling that $\mathcal{K}=\mathcal{K}_1\times\ldots\times\mathcal{K}_k$ is a Cartesian product of cones $\mathcal{K}_i$ for $i=1,\ldots,k$, let $F_i:\interior\mathcal{K}_i\rightarrow\mathbb{R}$ be a $\nu_i$-logarithmically homogeneous self-concordant barrier for $\mathcal{K}_i$. Then a $\nu$-logarithmically homogeneous self-concordant barrier for $\mathcal{K}$ is
\begin{equation*}
    F(s) = \sum_{i=1}^k F_i(s_i),
\end{equation*}
where $\nu=\sum_{i=1}^k \nu_i$ and we similarly split $s$ as $s_i=(s_1,\ldots,s_k)$ where $s\in\interior\mathcal{K}_i$ (and similarly for $z$). Let us denote $\Omega = \mathbb{R}^n\times\mathbb{R}^p\times\mathcal{K}_*\times\mathcal{K}\times\mathbb{R}_{+}\times\mathbb{R}_{+}$. For an initial point $\omega_0 \in \interior\Omega$, the central path is defined as the trajectory of points $\{\omega(\mu) : 0<\mu\leq1\}$ satisfying
\begin{equation}\label{eqn:central-path}
    \begin{aligned}
         L(\omega(\mu)) &= \mu L(\omega_0) \\
         z(\mu) + \mu \nabla F(s(\mu)) &= 0 \\
        \tau(\mu)\kappa(\mu) &= \mu \\
        (s(\mu), \kappa(\mu), z(\mu), \tau(\mu)) &\in \interior\mathcal{K} \times \mathbb{R}_{++} \times \interior\mathcal{K}_* \times \mathbb{R}_{++}.
    \end{aligned}
\end{equation}
As $\mu\downarrow0$, the central path converges to the solution of the original problem~\eqref{eqn:hsde}~\cite{nesterov2018lectures}. Primal-dual interior point methods track the central path to this solution by linearizing the central path equations to find suitable step directions.

\subsection{Initialization}\label{eqn:initialization}

Before running the primal-dual interior point algorithm, we first preprocess the problem using the same method as Hypatia by rescaling $c$, $A$, $b$, $G$ and $h$ to try make all values approximately the same magnitude. By default, we then initialize the interior point algorithm with initial points $\tau^0=\kappa^0=1$, $x^0=y^0=0$, and (where possible) choose $s^0$ and $z^0$ to satisfy
\begin{equation*}
    s^0 = z^0 = -\nabla F(s^0).
\end{equation*}
This is cheap to approximate (e.g., by solving a small nonlinear system of equations) or precompute for all cones we currently implement except the quantum key distribution cone $\mathcal{QKD}_{\mathcal{G}, \mathcal{Z}}$, which we instead initialize using the same method as~\cite{lorente2024quantum}.

\subsection{Stepping algorithm}

To follow the central path~\eqref{eqn:central-path} towards the solution, we implement two different stepping algorithms depending on whether $\mathcal{K}$ is a symmetric or nonsymmetric cone, i.e., the Nesterov-Todd algorithm~\cite{nesterov1997self,nesterov1998primal} if $\mathcal{K}$ is symmetric, and the Skajaa-Ye algorithm~\cite{skajaa2015homogeneous} if $\mathcal{K}$ is nonsymmetric. We outline these algorithms in more detail in the following subsections.

\subsubsection{Symmetric algorithm}

When $\mathcal{K}$ is a symmetric cone, i.e., a Cartesian product of the nonnegative orthant, positive semidefinite cone, and second-order cone, we use CVXOPT's~\cite{vandenberghe2010cvxopt} implementation of the Nesterov-Todd algorithm. For a given iterate $\omega\in\interior\Omega$, let us define a scaling point $w\in\mathbb{R}^q$ which is the unique point satisfying
\begin{equation*}
    \nabla^2 F(w) s = z,
\end{equation*}
a primal-dual scaling matrix $W$ which satisfies $\nabla^2 F(w)^{-1} = W^\top W$, and a scaled variable $\lambda=W^{-\top} s = Wz$. We refer to~\cite[Section 4]{vandenberghe2010cvxopt} for further details about these objects and how to efficiently compute them. Given these, we define the prediction direction $\Delta \omega_\text{p}$ as the solution to
\begin{equation}
    \begin{aligned}
        L(\Delta \omega_\text{p}) &= -L(\omega) \\
        \Delta z_\text{p} + \grad^2 F(w) \Delta s_\text{p} &= -z \\
        \tau \Delta \kappa_\text{p} + \kappa \Delta\tau_\text{p} &= -\tau\kappa,
    \end{aligned}
\end{equation}
and compute a centering parameter $\sigma=(1-\alpha_\text{p})^3$ where $\alpha_\text{p}$ is the largest step size in the direction $\Delta \omega_\text{p}$ to the boundary of the cone, i.e.,
\begin{equation*}
    \alpha_\text{p} = \max \{ \alpha \in [0, 1] : \omega + \alpha\Delta \omega_p \in \Omega \}.
\end{equation*}
See~\cite[Section 8]{vandenberghe2010cvxopt} for details about how this step size can be efficiently computed. Following this, we define the combined direction $\Delta \omega$ as the solution to 
\begin{equation}
    \begin{aligned}
        L(\Delta \omega) &= (1-\sigma) L(\omega) \\
        \Delta z + \grad^2 F(w) \Delta s &= W^\top (\lambda \diamond ( -\lambda \circ \lambda - (W^{-\top}\Delta s_\text{p}) \circ (W\Delta z_\text{p}) + \sigma \mu \bm{e} )) \\
        \tau \Delta \kappa + \kappa \Delta\tau &= -\tau\kappa - \Delta\tau_\text{p}\Delta\kappa_\text{p} + \sigma\mu,
    \end{aligned}
\end{equation}
where $\circ$ denotes the Jordan product on $\mathcal{K}$, $\bm{e}$ denotes the standard idempotent on $\mathcal{K}$, and $\diamond$ denotes the inverse of $\circ$ such that $u \circ (u\diamond v) = v$ for all $u\in\interior\mathcal{K}$ and $v\in\mathbb{R}^q$. See~\cite[Section 3 and Section 5.4]{vandenberghe2010cvxopt} for concrete details about these objects. Using this combined step, we update our iterates as follows
\begin{equation*}
    \omega^+ = \omega + 0.99 \alpha \Delta \omega,
\end{equation*}
where $\alpha$ is the largest step size in the direction $\Delta \omega$ to the boundary of the cone, i.e.,
\begin{equation*}
    \alpha = \max \{ \alpha \in [0, 1] : \omega + \alpha\Delta \omega \in \Omega \}.
\end{equation*}

\subsubsection{Nonsymmetric algorithm}\label{subsec:nonsym-alg}

When $\mathcal{K}$ is nonsymmetric, we instead use the combined stepping variant of the Skajaa-Ye algorithm as used by Hypatia~\cite{coey2023performance}. At a given interior point $\omega \in \interior\Omega$, let us define
\begin{equation*}
    \mu(\omega) = \frac{s^\top z + \tau\kappa}{\nu+1}.
\end{equation*}
Similar to the symmetric algorithm, we first define the prediction direction $\Delta \omega_\text{p}$ as the solution to
\begin{equation}
    \begin{aligned}
        L(\Delta \omega_\text{p}) &= -L(\omega) \\
        \Delta z_\text{p} + \mu(\omega) \grad^2 F(s) \Delta s_\text{p} &= -z \\
        \tau \Delta \kappa_\text{p} + \kappa \Delta\tau_\text{p} &= -\tau\kappa,
    \end{aligned}
\end{equation}
and its third order adjustment $\Delta \omega_\text{p}^{\text{toa}}$ as the solution to
\begin{equation}
    \begin{aligned}
        L(\Delta \omega_\text{p}^{\text{toa}}) &= 0 \\
        \Delta z_\text{p}^{\text{toa}} + \mu(\omega) \grad^2 F(s) \Delta s_\text{p}^{\text{toa}} &= \mu(\omega) \grad^2 F(s) \Delta s_\text{p} + \mu(\omega) \mathsf{D}^3 F(s)[s, \Delta s_\text{p}] \\
        \tau \Delta \kappa_\text{p}^{\text{toa}} + \kappa \Delta\tau_\text{p}^{\text{toa}} &= \Delta\tau_\text{p} \Delta \kappa_\text{p}.
    \end{aligned}
\end{equation}
Additionally, we define the centering direction $\Delta \omega_\text{c}$ as the solution to
\begin{equation}
    \begin{aligned}
        L(\Delta \omega_\text{c}) &= 0 \\
        \Delta z_\text{c} + \mu(\omega) \grad^2 F(s) \Delta s_\text{c} &= -z - \mu(\omega) \grad F(s) \\
        \tau \Delta \kappa_\text{c} + \kappa \Delta\tau_\text{c} &= -\tau\kappa + \mu(\omega),
    \end{aligned}
\end{equation}
and its third order adjustment $\Delta \omega_\text{c}^{\text{toa}}$ as the solution to
\begin{equation}
    \begin{aligned}
        L(\Delta \omega_\text{c}^{\text{toa}}) &= 0 \\
        \Delta z_\text{c}^{\text{toa}} + \mu(\omega) \grad^2 F(s) \Delta s_\text{c}^{\text{toa}} &= \mu(\omega) \mathsf{D}^3 F(s)[s, \Delta s_\text{c}] \\
        \tau \Delta \kappa_\text{c}^{\text{toa}} + \kappa \Delta\tau_\text{c}^{\text{toa}} &= 0.
    \end{aligned}
\end{equation}
Using these stepping directions, we update our iterates by taking the combined step
\begin{equation*}
    \omega^+ = \omega + \alpha (\Delta \omega_\text{p} + \alpha \Delta \omega_\text{p}^{\text{toa}}) + (1 - \alpha) (\Delta \omega_\text{c} + (1 - \alpha) \Delta \omega_\text{c}^{\text{toa}}),
\end{equation*}
where $\alpha$ is chosen using a backtracking line search to approximately satisfy
\begin{equation*}
    \alpha = \max\{ \alpha \in [0, 1] : \omega + \Delta \omega \in \mathcal{N}(0.99) \},
\end{equation*}
and we define the neighborhood $\mathcal{N}(\eta)$ for $0<\eta<1$ as
\begin{equation}\label{eqn:nonsym-proximity}
    \mathcal{N}(\eta) = \{ \omega \in \interior \Omega : \norm{ \grad^2 F_i(s_i)^{-1/2} ( z_i/\mu(\omega) + \grad F_i(s_i) ) } \leq \eta, \ \forall i=1,\ldots,k \}.
\end{equation}
From~\cite[Lemma 1]{coey2023performance}, we know that it suffices to check that $s\in\interior\mathcal{K}$ and the norm inequalities to ensure that $\omega\in\mathcal{N}(\eta)$, i.e., we do not need to explicitly check if $z\in\interior\mathcal{K}_*$.

\subsection{Solving the Newton system}\label{sec:newton}

For all algorithms and stepping directions, the main work is in solving linear systems of the form
\begin{equation}\label{eqn:newton}
    \begin{gathered}
        L(\Delta \omega) \coloneqq\begin{bmatrix}
            0 & A^\top & G^\top & c \\
            -A & 0 & 0 & b \\
            -G & 0 & 0 & h \\
            -c^\top & -b^\top & -h^\top & 0
        \end{bmatrix} \begin{bmatrix}
            \Delta x \\ \Delta y \\ \Delta z \\ \Delta\tau
        \end{bmatrix} - \begin{bmatrix}
            0 \\ 0 \\ \Delta s \\ \Delta\kappa
        \end{bmatrix} = \begin{bmatrix}
            r_x \\ r_y \\ r_z \\ r_\tau
        \end{bmatrix} \\
        \Delta z + H \Delta s = r_s \\
        \tau\Delta\kappa + \kappa\Delta\tau = r_\kappa,
    \end{gathered}
\end{equation}
where $H=\grad^2F(w)$ if we use the symmetric algorithm, and $H=\mu(\omega) \grad^2F(s)$ if we use the nonsymmetric algorithm. A suitable block elimination of this system of equations shows that the solution to this system of equations is
\begin{equation}
    \begin{aligned}
        \Delta \tau &= \frac{r_\tau + r_\kappa/\tau + c^\top \Delta x_{1} + b^\top \Delta y_{1} + h^\top \Delta z_{1}}{\kappa/\tau + c^\top \Delta x_{2} + b^\top \Delta y_{2} + h^\top \Delta z_{2}} \\
        \Delta x &= \Delta x_{1} - \Delta\tau \Delta x_{2} \\
        \Delta y &= \Delta y_{1} - \Delta\tau \Delta y_{2} \\
        \Delta z &= \Delta z_{1} - \Delta\tau \Delta z_{2} \\
        \Delta s &= -G\Delta x + \Delta\tau h - r_z \\
        \Delta \kappa &= (r_\kappa - \kappa\Delta\tau)/\tau,
    \end{aligned}
\end{equation}
where $(\Delta x_{1}, \Delta y_{1}, \Delta z_{1})$ and $(\Delta x_{2}, \Delta y_{2}, \Delta z_{2})$ are the solutions to the linear systems
\begin{equation}\label{eqn:newton-system-elim-a}
    \begin{bmatrix}
        0 & A^\top & G^\top \\
        -A & 0 & 0 \\
        -HG & 0 & \mathbb{I}
    \end{bmatrix} \begin{bmatrix}
        \Delta x_{i} \\ \Delta y_i \\ \Delta z_i
    \end{bmatrix} = \begin{bmatrix}
        b_{x,i} \\ b_{y,i} \\ b_{z,i}
    \end{bmatrix},
\end{equation}
where $(b_{x,1}, b_{y,1}, b_{z,1})=(r_x, r_y, Hr_z + r_s)$, and $(b_{x,2}, b_{y,2}, b_{z,2})=(c, b, Hh)$. To solve this $3\times3$ block subsystem, another block elimination gives us the normal equations
\begin{equation}\label{eqn:newton-system-elim-b}
    \begin{aligned}
        A (G^\top H G)^{-1} A^\top \Delta y_i &= b_{y,i} + A (G^\top H G)^{-1} (b_{x,i} - G^\top b_{z,i}) \\
        \Delta x_i &= (G^\top H G)^{-1} (b_{x,i} - G^\top b_{z,i} - A^\top \Delta y_i) \\
        \Delta z_i &= b_{z,i} + HG\Delta x_i,
    \end{aligned}
\end{equation}
i.e., the main work in solving the Newton system is in building and Cholesky factoring the matrices $G^\top H G$ and $A (G^\top H G)^{-1} A^\top$. If $G^\top H G$ is singular, then following~\cite[Section 10.1]{vandenberghe2010cvxopt}, we instead use the modified normal equations
\begin{equation}
    \begin{aligned}
        A (G^\top H G + A^\top A)^{-1} A^\top \Delta y_i &= b_{y,i} + A (G^\top H G + A^\top A)^{-1} (b_{x,i} - G^\top b_{z,i} - A^\top b_{y, i}) \\
        \Delta x_i &= (G^\top H G + A^\top A)^{-1} (b_{x,i} - G^\top b_{z,i} - A^\top (b_{y, i} + \Delta y_i)) \\
        \Delta z_i &= b_{z,i} + HG\Delta x_i.
    \end{aligned}
\end{equation}
When $G=-\mathbb{I}$, i.e., our conic constraints are of the form $x+h\in\mathcal{K}$, then the normal equations reduce to
\begin{equation}\label{eqn:newton-simplified}
    \begin{aligned}
        A H^{-1} A^\top \Delta y_i &= b_{y,i} + A H^{-1} (b_{x,i} + b_{z,i}) \\
        \Delta x_i &= H^{-1} (b_{x,i} - A^\top \Delta y_i + b_{z,i}) \\
        \Delta z_i &= A^\top \Delta y_i - b_{x,i},
    \end{aligned}
\end{equation}
i.e., now the main work is in building and Cholesky factoring the matrix $A H^{-1} A^\top$. If the residuals of the solution to the Newton system are large, then we perform iterative refinement to improve the quality of the step directions.

\subsection{Termination criteria}\label{sec:termination}

We use a similar stopping criteria as MOSEK~\cite{mosek}, which we outline below. First, we terminate the algorithm with an approximately primal-dual optimal solution $(x, y, z)/\tau$ if the current iterate $\omega$ has approximately zero duality gap, i.e., satisfies
\begin{equation*}
    \min\left( \frac{s^\top z}{\tau}, \abs{ c^\top x + b^\top y + h^\top z } \right) \leq \varepsilon_g \max\left(\tau, \frac{\min(\abs{c^\top x}, \abs{b^\top y + h^\top z})}{\tau} \right),
\end{equation*}
and is approximately primal-dual feasible, i.e., satisfies
\begin{align*}
    \norm{c + A^\top y/\tau + G^\top z/\tau}_\infty &\leq \varepsilon_f (1 + \norm{c}_\infty) \\ 
    \norm{b - Ax/\tau}_\infty &\leq \varepsilon_f (1 + \norm{b}_\infty) \\ 
    \norm{h - Gx/\tau - s/\tau}_\infty &\leq \varepsilon_f (1 + \norm{h}_\infty).
\end{align*}
We terminate the algorithm with an approximate certificate for primal infeasibility, i.e., a dual improving ray $(y, z)$, if $b^\top y + h^\top z < 0$ and
\begin{align*}
    \norm{A^\top y + G^\top z}_\infty \leq - \varepsilon_{i} (b^\top y + h^\top z).
\end{align*}
Similarly, we terminate the algorithm with an approximate certificate for dual infeasibiltiy, i.e., a primal improving ray $x$, if $c^\top x < 0$ and
\begin{align*}
    \max(\norm{Ax}_\infty, \norm{Gx + s}_\infty) \leq - \varepsilon_{i} (c^\top x).
\end{align*}
Finally, we claim the conic problem is ill-posed if $(x, y, z) \neq 0$ and
\begin{equation*}
    \max(\norm{A^\top y + G^\top z}_\infty, \norm{Ax}_\infty, \norm{Gx + s}_\infty) \leq \varepsilon_{ip} \max(\norm{x}_\infty, \norm{y}_\infty, \norm{z}_\infty).
\end{equation*}

\subsection{Avoiding inverse Hessian-vector products}\label{eqn:avoid-inv-hess}

When using the Skajaa-Ye variant of the interior-point algorithm described in Section~\ref{subsec:nonsym-alg}, when $G\neq-\mathbb{I}$, the only step which requires solving linear systems with $H$ is to compute the proximity measure~\eqref{eqn:nonsym-proximity}. All other steps only require solving linear systems with the matrix $G^\top H G$. When $G$ is a tall matrix, i.e., the dimension of $x$ is small while the dimension of $\mathcal{K}$ is large, $G^\top H G$ can be significantly smaller than $H$. When, additionally, we do not have efficient methods for computing inverse Hessian-vector products $H^{-1}v$, such as for the quantum relative entropy cone or noncommutative perspective cones, solving linear systems with $G^\top H G$ can be much easier than solving linear systems with $H$. Therefore, it is desirable to have a method which works directly in the subspace of $G$.

For simplicity, let us assume that $h=0$ (in Remark~\ref{rem:no-h}, we show that we can make this assumption without loss of generality). Assume $G$ has full column rank, and $\{ x\in\mathbb{R}^n : Gx \in \interior\mathcal{K} \}\neq \varnothing$. Rather than solving the original problem~\eqref{eqn:primal}, we consider solving the following equivalent conic program
\begin{gather}\label{eqn:simplified-primal}
    \begin{aligned}
        \minimize_{x\in\mathbb{R}^n} \quad & c^\top x \\
        \subjto \quad & b - Ax = 0 \\
                \quad & -x \in G^{-1}(\mathcal{K}) \coloneqq \{ x\in\mathbb{R}^n : Gx \in \mathcal{K} \}.
    \end{aligned}
\end{gather}
From~\cite[Proposition 2.3.3(i)]{nesterov1994interior}, we know that if $F:\mathbb{R}^q\rightarrow\mathbb{R}$ is a $\nu$-logarithmically homogeneous self-concordant barrier for $\mathcal{K}$, then $F'(x) = F(Gx)$ is a $\nu$-logarithmically homogeneous self-concordant barrier for $G^{-1}(\mathcal{K})$. Therefore, by working with this new cone and barrier, all steps of the Skajaa-Ye algorithm, including computing the proximity measure, now only require solving linear systems with the matrix $G^\top H G$, as desired.

To implement this strategy in practice, rather than requiring users to define a completely new cone oracle for $G^{-1}(\mathcal{K})$, we instead use the following observation. Let $\omega=(x, y, z, Gs', \tau, \kappa)$ be an interior point corresponding to the original problem~\eqref{eqn:primal} for some $s'\in\mathbb{R}^n$, and let $\omega'=(x', y', G^\top z, s', \tau', \kappa')$ be an interior point corresponding to the modified problem~\eqref{eqn:simplified-primal}. One can verify that the Newton directions for the Skajaa-Ye algorithm obtained at these points for each corresponding problem are equivalent, in the sense that 
\begin{equation*}
    (\Delta x, \Delta y, G^\top \Delta z, \Delta s, \Delta \tau, \Delta \kappa)=(\Delta x', \Delta y',  \Delta z', G \Delta s', \Delta \tau', \Delta \kappa').
\end{equation*}
Therefore, let us assume that we can easily obtain an initial point $\omega^0=(x^0, y^0, z^0, s^0, \tau^0, \kappa^0)$ satisfying $s^0\in\interior\mathcal{K}\cap \image G$ (we show in Remark~\ref{rem:s0} that this is always possible). Then to solve the modified conic program~\eqref{eqn:simplified-primal} using the Skajaa-Ye algorithm, we instead solve the original conic program~\eqref{eqn:primal} using the Skajaa-Ye algorithm initialized at $\omega^0$, and use the modified proximity measure
\begin{equation*}
    \mathcal{N}'(\eta) = \{ \omega\in \interior\tilde{\Omega} : \norm{ (G^\top\grad^2 F(s)G)^{-1/2} G^\top ( z/\mu(\omega) + \grad F(s) ) } \leq \eta\},
\end{equation*}
instead of the original proximity measure~\eqref{eqn:nonsym-proximity}, where $\tilde{\Omega}=\mathbb{R}^n\times\mathbb{R}^p\times \mathbb{R}^q \times\mathcal{K}\cap \image G\times\mathbb{R}_{+}\times\mathbb{R}_{+}$. A similar argument as~\cite[Lemma 1]{coey2023performance} (see, also, \cite[Lemma 15]{papp2017homogeneous}) shows that any point $\omega\in\mathcal{N}'(\eta)$ for $0<\eta<1$ satisfies $G^\top z\in\interior G^{-1}(\mathcal{K})_*$. Overall, this strategy allows us to indirectly solve the modified cone program~\eqref{eqn:simplified-primal} without having to significantly modify our algorithm or cone oracles, while still benefiting from no longer requiring inverse Hessian-vector products for $\mathcal{K}$.

\begin{rem}\label{rem:no-h}
    One of the assumptions we made in the above discussion was that the original problem~\eqref{eqn:primal} satisfied $h=0$. Here, we show we can always form an equivalent conic problem which satisfies $h'=0$. If $h\notin\image G$, then we can instead solve the equivalent conic program 
    \begin{equation*}
        \minimize_{x\in\mathbb{R}^n, \alpha\in\mathbb{R}} \quad c^\top x \quad \subjto \quad Ax = b, \quad \alpha=1, \quad \alpha h - Gx \in \mathcal{K}.
    \end{equation*}    
    Note if $G$ is full column rank, then the matrix $G'=(G, h)$ is also full column rank. If $h\in\image G$, then we can eliminate $h$ by performing a change of variables from $x$ to $z=x-G^+h$, where $G^+$ is the pseudo-inverse of $G$, i.e.,
    \begin{equation*}
        \minimize_{z\in\mathbb{R}^n} \quad c^\top z + c^\top(G^+h) \quad \subjto \quad Az = b - A(G^+h), \quad -Gz \in \mathcal{K}.
    \end{equation*}
\end{rem}

\begin{rem}\label{rem:s0}
    A second assumption we made was that we could easily obtain an initial point $s^0\in\interior\mathcal{K}\cap \image G$ for the original problem \eqref{eqn:primal}. Here, we show that we always form an equivalent conic problem such that $s^0\in\image G'$ for any given $s^0\in\interior\mathcal{K}$. If $s^0$ is already in the image of $G$, then this is trivial. If $s^0\notin\image G$, then we instead solve the equivalent cone program
    \begin{equation*}
        \minimize_{x\in\mathbb{R}^n, \beta\in\mathbb{R}} \quad c^\top x \quad \subjto \quad Ax = b, \quad \beta=0, \quad \beta s^0 - Gx \in \mathcal{K}.
    \end{equation*}    
    This trivially guarantees that $s^0$ is in the image of the matrix $G'=(G, s^0)$. Additionally, $G'$ is full column rank if $G$ is full column rank.    
\end{rem}


\section{Cone oracles}\label{sec:cones}

To perform the primal-dual interior point algorithm discussed in Section~\ref{sec:pdipm}, we require efficient oracles for the gradient, Hessian-vector product, inverse Hessian-vector product, and third order derivatives of the barrier function $F$ of the cone $\mathcal{K}$. QICS provides implementations of these oracles for a variety of cones which commonly arise in quantum information theory, which we have listed in Section~\ref{subsec:native-cones}. In the following, we provide implementation details for the family of quantum entropy and noncommutative perspective cones which we support. We omit details about third order derivatives as these are relatively straightforward to derive. In Section~\ref{subsec:sparsity}, we also discuss how we exploit sparsity in the problem data, which in practice can significantly improve the efficiency of the interior-point algorithm.

\subsection{Preliminaries}\label{sec:quantentr}


Consider a twice continuously differentiable scalar function $g:\domain g\subseteq\mathbb{R} \rightarrow \mathbb{R}$. We can extend this function to act on Hermitian matrices $X$ with eigendecomposition $X=\sum_{i=1}^n \lambda_i u_iu_i^\dag$ such that $\lambda_i\in\domain g$ for all $i=1,\ldots,n$ by defining
\begin{equation*}
    g(X) \coloneqq \sum_{i=1}^n g(\lambda_i) u_i u_i^\dag.
\end{equation*}
The derivatives of this function are~\cite[Theorem 3.33]{hiai2014introduction}
\begin{subequations}\label{eqn:deriv-spectral}
    \begin{align}
        \mathsf{D} g(X)[V] &= U[g^{[1]}(\Lambda) \odot (U^\dag V U)]U^\dag\\
        \mathsf{D}^2 g(X)[V,W] &= U \biggl[ \sum_{k=1}^n g^{[2]}_k(\Lambda) \odot \biggl([ U^\dag VU]_k[ U^\dag WU]_k^\dag + [U^\dag WU]_k[U^\dag VU]_k^\dag \biggr) \biggr]U^\dag,
    \end{align}
\end{subequations}
where $\Lambda=\diag(\lambda_1,\ldots,\lambda_n)$, $g^{[1]}(\Lambda)$ is the first divided differences matrix whose $(i,j)$-th entry is $g^{[1]}(\lambda_i, \lambda_j)$, where
\begin{align*}
    g^{[1]}(\lambda, \mu) &= \frac{g(\lambda) - g(\mu)}{\lambda - \mu}, \quad \textrm{if }\lambda\neq\mu\\
    g^{[1]}(\lambda, \lambda) &= g'(\lambda),
\end{align*} 
and $g^{[2]}_k(\Lambda)$ for $k=1,\ldots,n$ are the second divided difference matrices whose $(i,j)$-th entries are $g^{[2]}(\lambda_i, \lambda_j, \lambda_k)$, where
\begin{align*}
    g^{[2]}(\lambda, \mu, \sigma) &= \frac{g^{[1]}(\lambda, \sigma) - g^{[1]}(\mu, \sigma)}{\lambda - \mu}, \quad \textrm{if }\lambda\neq\mu \\
    g^{[2]}(\lambda, \lambda, \lambda) &= \frac{1}{2}g''(\lambda),
\end{align*}
and all other scenarios are given by defining $g^{[2]}(\lambda, \mu, \sigma)$ symmetrically in each of its arguments. Notably, if $\mathsf{D}g(X)$ is invertible, then $g^{[1]}(\Lambda)$ is element-wise nonzero, and we can easily take the inverse of the first derivative map by exchanging the elementwise multiplication for an elementwise division, i.e.,
\begin{align}\label{eqn:deriv-spectral-inverse}
    \mathsf{D} g(X)^{-1} [V] &= U[(U^\dag V U) \oslash g^{[1]}(\Lambda)]U^\dag.
\end{align}
We also recall that the derivatives of the log determinant $X\in\mathbb{H}^n_{++}\mapsto\log\det(X)$, which is the barrier function for the positive semidefinite cone $\mathbb{H}^n_+$, are
\begin{subequations}\label{eqn:logdet-derivatives}
    \begin{align}
        \mathsf{D} \log\det(X) &= X^{-1} \\
        \mathsf{D}^2 \log\det(X)[V] &= -X^{-1}VX^{-1}.
    \end{align}
\end{subequations}
Finally, we recall the chain rule, which we state below.
\begin{lem}[{\cite[Theorem 3.4]{higham2008functions}}]\label{lem:chain-rule}
    Consider the functions $g:\mathbb{V}\rightarrow\mathbb{V}'$ and $f:\mathbb{V}'\rightarrow\mathbb{V}''$ which are differentiable at $x$ and $g(x)$, respectively. Then for all $v\in\mathbb{V}$, we have
    \begin{equation*}
        \mathsf{D}(f \circ g)(x)[v] = \mathsf{D}f(g(x))[\mathsf{D}g(x)[v]].
    \end{equation*}
\end{lem}

\subsubsection{Generalized barrier function}

Consider the function $F:\mathbb{R}\times\mathbb{H}^n_{++}\times\mathbb{H}^m_{++}\rightarrow\mathbb{R}$ defined by
\begin{equation}\label{eqn:general-barrier}
    F(t, X, Y) = -\log(t - f(X, Y)) - \log\det(X) - \log\det(Y),
\end{equation}
for some convex function $f:\mathbb{H}^n_{++}\times\mathbb{H}^m_{++}\rightarrow\mathbb{R}$. This is the form of the barrier function for most of the cones we discuss in the following subsections. A straightforward application of the chain rule shows that the gradients of $F$ are given by
\begin{gather}\label{eqn:t-barrier-grad}
    \begin{aligned}
        \grad F(t, X, Y) = \begin{bmatrix}
            -z^{-1} \\
            z^{-1} \grad_X f(X, Y) - \grad \log\det(X) \\
            z^{-1} \grad_Y f(X, Y) - \grad \log\det(Y)
        \end{bmatrix},
    \end{aligned}
\end{gather}
where $z = t - f(X, Y)$. 
Similarly, one can show that the Hessian matrix of $F$ has the form
\begin{equation}\label{eqn:t-barrier-hessian}
    \grad^2F(t, X, Y)
    =
    \begin{bmatrix} z^{-1} \\ -z^{-1}\grad_X f (X,Y) \\ -z^{-1}\grad_Y f (X,Y) \end{bmatrix}\begin{bmatrix} z^{-1} \\ -z^{-1}\grad_X f (X,Y) \\ -z^{-1}\grad_Y f (X,Y) \end{bmatrix}^\top + \begin{bmatrix}
        0 & 0 & 0 \\
        0 & M_{11} & M_{12} \\
        0 & M_{21} & M_{22}
    \end{bmatrix},
\end{equation}
where $M$ is the $2\times2$ block matrix
\begin{equation}\label{eqn:general-barrier-system}
    M = z^{-1} \grad^2 f(X, Y) - \begin{bmatrix}
        \grad^2 \log\det(X) & 0 \\
        0 & \grad^2 \log\det(Y)
    \end{bmatrix}.
\end{equation}
Therefore, to solve the linear system of equations
\begin{equation*}
    \grad^2F(t, X, Y) \begin{bmatrix}
        t \\
        \vect(\Delta X) \\
        \vect(\Delta Y) \\
    \end{bmatrix} = \begin{bmatrix}
        r_t \\
        \vect(R_x) \\
        \vect(R_y) \\
    \end{bmatrix},
\end{equation*}
a suitable block elimination shows that 
\begin{subequations}\label{eqn:t-barrier-elim}
\begin{gather}
    M \begin{bmatrix}
        \vect(\Delta X) \\ \vect(\Delta Y)
    \end{bmatrix} = \begin{bmatrix}
        \vect(R_x) + z^{-1} \grad_X f(X, Y) \\ \vect(R_y) + z^{-1} \grad_Y f(X, Y)
    \end{bmatrix} \label{eqn:t-barrier-elim-a} \\
    \Delta t = z^2 r_t + \mathsf{D}_X f(X, Y)[\Delta X] + \mathsf{D}_Y f(X, Y)[\Delta Y].
\end{gather}
\end{subequations}
That is, the main work in solving linear systems with the Hessian matrix of the generalized barrier function~\eqref{eqn:general-barrier} is in solving linear systems with the $2\times2$ block matrix~\eqref{eqn:general-barrier-system}. 

We now outline some specific details about the implementations of the family of quantum entropy cones.

\subsection{Quantum entropy}
Recall that the (homogenized) quantum entropy cone is defined as
\begin{equation*}
    \mathcal{QE}_n \coloneqq \closure\{ (t, u, X) \in \mathbb{R}\times\mathbb{R}_{++}\times\mathbb{H}^n_{++} : t \geq -uS(u^{-1}X) \}.
\end{equation*}
In~\cite{coey2023conic}, it was shown that the function
\begin{equation}\label{eqn:qe-barrier}
    (t, u, X) \in \mathbb{R}\times\mathbb{R}_{++}\times\mathbb{H}^n_{++} \mapsto -\log(t + uS(u^{-1}X)) - \log(u) - \log\det(X),
\end{equation}
is a $(2+n)$-logarithmically homogeneous self-concordant barrier for the quantum entropy cone. To obtain the derivatives of this barrier function, we first we use~\eqref{eqn:deriv-spectral} to obtain the first and second derivatives of the quantum entropy
\begin{subequations}\label{eqn:qe-derivatives}
    \begin{align}
        \mathsf{D} S(X) &= -\log(X) - \mathbb{I}, \\
        \mathsf{D}^2 S(X)[V] &= -U_x [\log^{[1]}(\Lambda_x) \odot (U^\dag_x V U_x)] U_x^\dag,
    \end{align}
\end{subequations}
where $X$ has eigendecomposition $X=U_x\Lambda_xU_x^\dag$. The derivatives of the homogenized quantum entropy and barrier function~\eqref{eqn:qe-barrier} are a straightforward modification of these results (see, also,~\cite{coey2023conic}). Notably, the main linear system of equations~\eqref{eqn:t-barrier-elim-a} we need to solve can be expressed as
\begin{equation*}
    \begin{bmatrix}
        z^{-1}u^{-2}\tr[X] + u^{-2} & -z^{-1}u^{-1}\vect(\mathbb{I})^\top \\
        -z^{-1}u^{-1}\vect(\mathbb{I}) & \grad^2 h(X)
    \end{bmatrix} \begin{bmatrix}
        \Delta u \\
        \vect(\Delta X)
    \end{bmatrix} = \begin{bmatrix}
        r_u \\
        \vect(R_x)
    \end{bmatrix},
\end{equation*}
where $h(x)=z^{-1}x\log(x) - \log(x)$. By recognizing that we can easily take the inverse of $\grad^2 h(X)$ by using~\eqref{eqn:deriv-spectral-inverse}, we can show that this system has the explicit solution
\begin{align*}
    \Delta u &= \frac{r_u + z^{-1}\tr[\Delta X_1]}{u^{-2}(1 + z^{-1}\tr[X]) - z^{-2}\tr[\Delta X_2]} \\
    \Delta X &= \Delta X_1 + uz^{-1}\Delta X_2,
\end{align*}
where
\begin{align*}
    \Delta X_1 &= U_x [ (U^\dag_x R_x U_x) \oslash h'^{[1]}(\Lambda_x) ] U_x^\dag \\
    \Delta X_2 &= U_x [ \mathbb{I} \oslash h'^{[1]}(\Lambda_x) ] U_x^\dag.
\end{align*}

\subsection{Quantum relative entropy} \label{subsec:cone-qre}
Recall that the quantum relative cone is defined as 
\begin{equation*}
    \mathcal{QRE}_n \coloneqq \{ (t, X, Y) \in \mathbb{R}\times\mathbb{H}^n_{+}\times\mathbb{H}^n_{+} : t \geq S\divx{X}{Y}, \ \ker(Y)\subseteq\ker(X) \}.
\end{equation*}
In~\cite{fawzi2023optimal}, it was shown that the function
\begin{equation}\label{eqn:qre-barrier}
    (t, X, Y)\in\mathbb{R}\times\mathbb{H}^n_{++}\times\mathbb{H}^n_{++} \mapsto -\log(t - S\divx{X}{Y}) - \log\det(X) - \log\det(Y),
\end{equation}
is a $(1+2n)$-logarithmically homogeneous self-concordant barrier for the quantum relative entropy cone. To obtain the derivatives of this barrier function, using~\eqref{eqn:deriv-spectral} we can show that the first derivatives of quantum relative entropy are
\begin{align*}
    \mathsf{D}_X S\divx{X}{Y} &= \log(X) - \log(Y) + \mathbb{I} \\
    \mathsf{D}_Y S\divx{X}{Y} &= -U_y [\log^{[1]}(\Lambda_Y) \odot (U^\dag_y X U_y)] U_y^\dag,
\end{align*}
where $X$ and $Y$ have eigendecompositions $X=U_x\Lambda_xU_x^\dag$ and $Y=U_y\Lambda_yU_y^\dag$, respectively. Similarly, we can show that the second derivatives of quantum relative entropy are
\begin{align*}
    \mathsf{D}^2_{XX} S\divx{X}{Y}[V_x] &= U_x [\log^{[1]}(\Lambda_x) \odot (U^\dag_x V_x U_x)] U_x^\dag \\
    \mathsf{D}^2_{XY} S\divx{X}{Y}[V_y] &= -U_y [\log^{[1]}(\Lambda_y) \odot (U^\dag_y V_y U_y)] U_y^\dag \\
    \mathsf{D}^2_{YX} S\divx{X}{Y}[V_x] &= -U_y [\log^{[1]}(\Lambda_y) \odot (U^\dag_y V_x U_y)] U_y^\dag \\
    \mathsf{D}^2_{YY} S\divx{X}{Y}[V_y] &= -U_y \biggl[\sum_{k=1}^n \log^{[2]}_k(\Lambda_y) \odot \biggl([U_y^\dag X U_y]_k[U_y^\dag V_y U_y]_k^\dag + [U_y^\dag V_y U_y]_k[U_y^\dag X U_y]_k^\dag \biggr) \biggr] U_y^\dag.
\end{align*}
Now to solve linear systems with the Hessian matrix of the barrier function, the main linear system of equations~\eqref{eqn:t-barrier-elim-a} we need to solve can be represented as
\begin{align*}
    \begin{bmatrix}
       \mathcal{U}_x & 0 \\ 0 & \mathcal{U}_y
    \end{bmatrix} \begin{bmatrix}
        \mathcal{D}_{xx} & -\mathcal{U}_{x}^\top \mathcal{U}_{y} \mathcal{D}_{xy} \\ -\mathcal{D}_{xy}\mathcal{U}_{y}^\top\mathcal{U}_{x} & \mathcal{S}_{yy}
    \end{bmatrix} \begin{bmatrix}
       \mathcal{U}_x & 0 \\ 0 & \mathcal{U}_y
    \end{bmatrix}^\top \begin{bmatrix}
        \vect(\Delta X) \\ \vect(\Delta Y)
    \end{bmatrix} = \begin{bmatrix}
        \vect(R_x) \\ \vect(R_y)
    \end{bmatrix},
\end{align*}
where $\mathcal{D}_{xx}$ and $\mathcal{D}_{xy}$ are the matrix representations of the linear maps
\begin{align*}
    V &\mapsto V \odot (z^{-1}\log^{[1]}(\Lambda_x) - \inv^{[1]}(\Lambda_x))) \\
    V &\mapsto V \odot (z^{-1}\log^{[1]}(\Lambda_y)),
\end{align*}
respectively, where $\inv(x)=1/x$. These linear maps should be interpreted as diagonal matrices which can be easily inverted following~\eqref{eqn:deriv-spectral-inverse}. Additionally, $\mathcal{U}_{x}$ and $\mathcal{U}_{y}$ are the matrix representations of the linear maps
\begin{align*}
    V &\mapsto U_x V U_x^\dag \\
    V &\mapsto U_y V U_y^\dag,
\end{align*}
respectively, and should be interpreted as orthogonal matrices such that $\mathcal{U}_x^{-1}=\mathcal{U}_x^\top$ and $\mathcal{U}_y^{-1}=\mathcal{U}_y^\top$. Finally, $\mathcal{S}_{yy}$ is a sparse matrix (see, e.g.,~\cite[Section 4.1.1]{faybusovich2020self}) representing the linear map
\begin{equation*}
    V \mapsto -z^{-1} \sum_{k=1}^n  ([U_y^\dag X U_y]_k V_k^\dag + V_k[U_y^\dag X U_y]_k^\dag) \odot \log^{[2]}_k(\Lambda_y) - V \odot \inv^{[1]}(\Lambda_y).
\end{equation*}
Given these, a suitable block elimination shows that the solutions $\Delta X$ and $\Delta Y$ are given by
\begin{align*}
    \vect(U_y^\dag \Delta Y U_y) &= \mathcal{C}^{-1} (\vect(U_y^\dag R_y U_y) - \mathcal{D}_{xy} \mathcal{U}_y^\top \mathcal{U}_x \mathcal{D}_{xx}^{-1} \vect(U_x^\dag R_x U_x) ) \\
    \vect(U_x^\dag \Delta X U_x) &= -\mathcal{D}_{xx}^{-1} (\vect(U_x^\dag R_x U_x) + \mathcal{U}_x^\top \mathcal{U}_y \mathcal{D}_{xy} \vect(\Delta Y)),
\end{align*}
where $\mathcal{C}$ is the Schur complement matrix
\begin{align*}
    \mathcal{C} = \mathcal{S}_{yy} - \mathcal{D}_{xy} \mathcal{U}_y^\top \mathcal{U}_x \mathcal{D}_{xx}^{-1} \mathcal{U}_x^\top \mathcal{U}_y \mathcal{D}_{xy}.
\end{align*}
Compared to directly constructing and Cholesky factoring the whole $2\times2$ block matrix, as done in Hypatia~\cite{coey2023performance}, this block elimination method only requires us to construct and Cholesky factor a matrix of half the size, and therefore requires around $8$ times fewer floating point operations.

\subsection{Quantum conditional entropy} \label{subsec:cone-qce}

Consider a system composed of $r$ subsystems, where the $i$-th subsystem has dimension $n_i$, and let $n=\Pi_i n_i$. Recall that the quantum conditional cone was defined as
\begin{equation}
    \mathcal{QCE}_{\{n_i\}, j} \coloneqq \closure\{ (t, X) \in \mathbb{R}\times\mathbb{H}^{\Pi_in_i}_{++} : t \geq -S(X) + S(\tr_i(X)) \}.
\end{equation}
In~\cite{he2024exploiting}, it was shown that the function
\begin{equation}\label{eqn:qce-barrier}
    (t, X)\in\mathbb{R}\times\mathbb{H}^n_{++} \mapsto -\log(t + S(X) - S(\tr_i(X))) - \log\det(X),
\end{equation}
is a $(1+n)$-logarithmically homogeneous self-concordant barrier for the quantum conditional entropy cone. The derivatives of this barrier function can be found using the chain rule together with the derivatives of quantum entropy~\eqref{eqn:qe-derivatives}. We use the method described in~\cite[Example 5.7]{he2024exploiting} to solve linear systems with the Hessian matrix of this barrier function, i.e., we can represent the main linear system~\eqref{eqn:t-barrier-elim-a} as
\begin{equation*}
    (\grad^2h(X) + z^{-1}\bm{\tr_i}^\top \grad^2 S(\tr_i(X)) \bm{\tr_i}) \vect(\Delta X) = \vect(R_x),
\end{equation*}
where $h(x)=z^{-1}x\log(x) - \log(x)$ and $\bm{\tr_i}$ is the matrix representation of $\tr_i$. As $\grad^2 h(X)$ is easy to solve linear systems with, and $\grad^2 S(\tr_i(X))$ is a small matrix, we can efficiently solve the required system of equations by using the matrix inversion lemma (see, e.g.,~\cite[Appendix A]{he2024exploiting}).

\subsection{Quantum key distribution} \label{subsec:cone-qkd}

Consider a positive linear map $\mathcal{G}:\mathbb{H}^n\rightarrow\mathbb{H}^{mr}$, and let $\mathcal{Z}:\mathbb{H}^{mr}\rightarrow\mathbb{H}^{mr}$ represent the pinching map that maps off-diagonal blocks of an $r\times r$ block matrix to zero. Recall that the quantum key distribution cone is defined as
\begin{equation}
    \mathcal{QKD}_{\mathcal{G},\mathcal{Z}} \coloneqq \closure \{ (t, X) \in \mathbb{R}\times\mathbb{H}^n_{++} : t \geq -S(\mathcal{G}(X)) + S(\mathcal{Z}(\mathcal{G}(X))) \},
\end{equation}
 In~\cite{he2024exploiting}, it was shown that the function 
\begin{equation}\label{eqn:qkd-barrier}
    (t, X) \in \mathbb{R}\times\mathbb{H}^n_{++} \mapsto -\log(t + S(\mathcal{G}(X)) - S(\mathcal{Z}(\mathcal{G}(X))) ) - \log\det(X),
\end{equation}
is a $(1+n)$-logarithmically homogeneous self-concordant barrier for the quantum conditional entropy cone. Like the conditional entropy cone, the derivatives of this barrier function can be found using the chain rule together with the derivatives of quantum entropy~\eqref{eqn:qe-derivatives}. To solve linear systems with the Hessian matrix of the barrier function we use the method described in~\cite[Section 6.1]{he2024exploiting}, i.e., we can exploit the fact that $\mathcal{Z}$ maps to block diagonal matrices to construct the required matrix~\eqref{eqn:general-barrier-system} more efficiently, which we then Cholesky factor to solve linear equations with. In the special case when $\mathcal{G}$ is the identity map, i.e., $\mathcal{G}(X)=X$, then the main linear system~\eqref{eqn:t-barrier-elim-a} reduces to
\begin{equation*}
    \biggl( \grad^2h(X) + z^{-1} \sum_{i=1}^r \bm{\mathcal{Z}_i}^\dag \grad^2 S(\mathcal{Z}_i(X)) \bm{\mathcal{Z}_i} \biggr) \vect(\Delta X) = \vect(R_x),
\end{equation*}
where $h(x)=z^{-1}x\log(x) - \log(x)$, and for each $i=1,\ldots,r$, $\mathcal{Z}_i:\mathbb{H}^{mr}\rightarrow\mathbb{H}^m$ is the linear map that outputs the $i$-th diagonal block of the $r\times r$ block matrix $X$ and $\bm{\mathcal{Z}_i}$ is its matrix representation. In this case, like the quantum conditional entropy cone, we can efficiently solve the required linear system of equations by using the matrix inversion lemma (see, e.g.,~\cite[Appendix A]{he2024exploiting}).

\subsection{Noncommutative perspective cones} \label{subsec:cone-op}

Consider an operator convex function $g:(0, \infty)\rightarrow\mathbb{R}$, meaning that for all $X,Y\in\mathbb{H}^n_{++}$, $\lambda\in[0, 1]$, and integers $n$, we have
\begin{equation*}
    g(\lambda X + (1-\lambda)Y) \preceq \lambda g(X) + (1-\lambda) g(Y).    
\end{equation*}
Recall that the epigraph of the noncommutative perspective of $g$ is defined as
\begin{equation}
    \mathcal{OPE}_n^g \coloneqq \closure\{ (T, X, Y) \in \mathbb{H}^n\times\mathbb{H}^n_{++}\times\mathbb{H}^n_{++} : T \succeq P_g(X, Y) \}.
\end{equation}
In~\cite{fawzi2023optimal}, it was shown that the function $F:\mathbb{H}^n\times\mathbb{H}^n_{++}\times\mathbb{H}^n_{++}\rightarrow\mathbb{R}$ defined by
\begin{equation}\label{eqn:op-barrier}
    F(T, X, Y) = -\log\det(T - P_g(X, Y)) - \log\det(X) - \log\det(Y),
\end{equation}
is a $3n$-logarithmically homogeneous self-concordant barrier for the epigraph of the noncommutative perspective. To derive the gradient of this barrier function, we can apply the chain rule introduced in Lemma~\ref{lem:chain-rule} together with the derivatives of the log determinant~\eqref{eqn:logdet-derivatives} to show that
\begin{align}
    \grad F(T, X, Y) = \begin{bmatrix}
        -\grad\log\det(Z)\\
        \grad\log\det(Z)^\top\grad_X P_g(X, Y) - \grad\log\det(X) \\
        \grad\log\det(Z)^\top\grad_Y P_g(X, Y) - \grad\log\det(Y)
    \end{bmatrix}
\end{align}
where $Z = T - f(X, Y)$.
Note the similarities between this gradient and that derived for the scalar case~\eqref{eqn:t-barrier-grad}. Using the identity $P_g(X,Y)=P_{\hat{g}}(Y,X)$ for all $X,Y\in\mathbb{H}^n_{++}$~\cite[Lemma 2.1]{hiai2017different}, where $\hat{g}(x)=xg(1/x)$ denotes the transpose of $g$, we can show using the chain rule that the required gradients of the noncommutative perspective are given by
\begin{align*}
    \mathsf{D}_X P_g(X, Y)[V_x] &= Y^{1/2} \mathsf{D}\hat{g}(\tilde{X})[\tilde{V}_x] Y^{1/2} \\
    \mathsf{D}_Y P_g(X, Y)[V_y] &= X^{1/2} \mathsf{D} g(\tilde{Y})[\tilde{V}_y] X^{1/2},
\end{align*}
where we denote $\tilde{X}=Y^{-1/2}XY^{-1/2}$, $\tilde{Y}=X^{-1/2}YX^{-1/2}$, $\tilde{V}_x=Y^{-1/2}V_xY^{-1/2}$ and $\tilde{V}_y=X^{-1/2}V_yX^{-1/2}$. Similarly, we can show that the Hessian matrix of $F$ is given by
\begin{equation*}
    \grad^2F(T, X, Y)
    =
    \begin{bmatrix}
        \mathcal{Z}^{-1} \\
        -\grad_X P_g(X,Y)^\top \mathcal{Z}^{-1} \\
        -\grad_Y P_g(X,Y)^\top \mathcal{Z}^{-1} \\
    \end{bmatrix}\begin{bmatrix}
        \mathcal{Z}^{-1} \\
        -\grad_X P_g(X,Y)^\top \mathcal{Z}^{-1} \\
        -\grad_Y P_g(X,Y)^\top \mathcal{Z}^{-1} \\
    \end{bmatrix}^\top
    + \begin{bmatrix}
        0 & 0 & 0 \\
        0 & M_{11} & M_{12} \\
        0 & M_{21} & M_{22}
    \end{bmatrix},
\end{equation*}
where $\mathcal{Z}^{-1}$ is the matrix representation of the linear map $V \mapsto L^{-\top} V L^{-1}$ where $L$ is the Cholesky factorization of $Z=LL^\top$, and $M$ is the $2\times2$ block matrix
\begin{equation}
    M = z^{-1} \grad^2 P_g(X, Y)^*[Z^{-1}] - \begin{bmatrix}
        \grad^2\log\det(X) & 0 \\
        0 & \grad^2\log\det(Y)
    \end{bmatrix},
\end{equation}
where $\grad^2P_g(X, Y)^*[Z^{-1}]$ is the matrix representation of the bilinear map
\begin{equation*}
    ((V_x, V_y), (W_x, W_y)) \mapsto \inp{Z^{-1}}{\mathsf{D}^2 P_g(X, Y)[(V_x, V_y), (W_x, W_y)]}.
\end{equation*}
Again, we remark on the similarities between this Hessian and the Hessian derived for the scalar case~\eqref{eqn:t-barrier-hessian}. Using the chain rule, we can show that the second derivatives of the noncommutative perspective are
\begin{align*}
    \mathsf{D}_{XX}^2 P_g(X, Y)[V_x, W_x] &= Y^{1/2} \mathsf{D}^2\hat{g}(\tilde{X})[\tilde{V}_x, \tilde{W}_x] Y^{1/2} \\
    \mathsf{D}_{YY}^2 P_g(X, Y)[V_y, W_y] &= X^{1/2} \mathsf{D}^2 g(\tilde{Y})[\tilde{V}_y, \tilde{W}_y] X^{1/2},
\end{align*}
where $\tilde{W}_x=Y^{-1/2}W_xY^{-1/2}$ and $\tilde{W}_y=X^{-1/2}W_yX^{-1/2}$. The second derivative for the cross-term $\mathsf{D}^2_{XY}P_g(X, Y)[V_x, V_y]$ is not as straightforward to obtain using a simple application of the chain rule. However, we show how this derivative can be obtained in closed form in the following lemma.
\begin{lem}
    Consider a twice continuously differentiable function $g:\domain g\rightarrow\mathbb{R}$. Let $h(x)=xg(x)$, and let $P_g:\mathbb{H}^n_{++}\times\mathbb{H}^n$ be the noncommutative perspective of $g$.
    \begin{equation*}
        \mathsf{D}^2_{XY}P_g(X, Y)[V_x, V_y] = X^{1/2} ( \tilde{V}_x \mathsf{D}g(\tilde{Y})[\tilde{V}_y] + \mathsf{D}g(\tilde{Y})[\tilde{V}_y] \tilde{V}_x - \mathsf{D}^2 h(\tilde{Y})[\tilde{V}_x, \tilde{V}_y]) X^{1/2},
    \end{equation*}
    where $\tilde{Y}=X^{-1/2}YX^{-1/2}$, $\tilde{V}_x=X^{-1/2}V_xX^{-1/2}$ and $\tilde{V}_y=X^{-1/2}V_yX^{-1/2}$.
\end{lem}
\begin{rem}
    In the scalar case, i.e., $n=1$, this second derivative simplifies to
    \begin{equation*}
        \frac{\partial^2 P_g}{\partial x \partial y} (x, y)[v_x, v_y] = -\frac{y}{x^2} g''\biggl(\frac{y}{x} \biggr) v_xv_y,
    \end{equation*}
    which is precisely what we expect when we directly take the derivative of the scalar perspective $(x,y)\mapsto x g(y/x)$ of $g$.
\end{rem}
\begin{proof}
    We will first consider the case when $g(x)=x^p$ for some $p\in\mathbb{N}$, then extend this result using a continuity argument. Using the chain rule together with~\cite[Example 3.32]{hiai2014introduction}, we obtain
    \begin{align*}
        \mathsf{D}_Y P_g(X, Y)[V_y] &= X^{1/2} \biggl( \sum_{k=1}^p \tilde{Y}^{k-1} \tilde{V_y} \tilde{Y}^{p-k} \biggr) X^{1/2} \\
        &= X \biggl( \sum_{k=1}^p X^{-1} (YX^{-1})^{k-1} V_y (X^{-1} Y)^{p-k} X^{-1} \biggr) X.
    \end{align*}
    Now using $\mathsf{D} f(X)[V]=-X^{-1}VX^{-1}$ when $f(X)=X^{-1}$ together with the product rule, we can show that
    \begin{align*}
        \mathsf{D}_{XY}^2 P_g(X, Y)[V_x, V_y] &= V_x X^{-1/2} \biggl( \sum_{k=1}^p \tilde{Y}^{k-1} \tilde{V_y} \tilde{Y}^{p-k} \biggr) X^{1/2} + X^{1/2} \biggl( \sum_{k=1}^p \tilde{Y}^{k-1} \tilde{V_y} \tilde{Y}^{p-k} \biggr) X^{-1/2} V_x \\
        &\hphantom{=} - X^{1/2} \biggl( \sum_{i+j+k=p-1} \tilde{Y}^i \tilde{V_x} \tilde{Y}^j \tilde{V_y} \tilde{Y}^k + \sum_{i+j+k=p-1} \tilde{Y}^i \tilde{V_y} \tilde{Y}^j \tilde{V_x} \tilde{Y}^k \biggr) X^{1/2} \\
        &= X^{1/2} ( \tilde{V}_x \grad g(\tilde{Y})[\tilde{V}_y] + \grad g(\tilde{Y})[\tilde{V}_y] \tilde{V}_x - \grad^2 h(\tilde{Y})[\tilde{V}_x, \tilde{V}_y]) X^{1/2},
    \end{align*}
    where the last equality comes by comparing to~\cite[Example 3.32]{hiai2014introduction}. This result extends to any polynomial $g$ due to linearity, and we can further extend this result to all twice continuously differentiable functions $g$ by using a similar continuity argument as~\cite[Theorem V.3.3]{bhatia2013matrix}.
\end{proof}
Finally, to solve linear systems of equations with the Hessian matrix of $F$, i.e.,
\begin{equation*}
    \grad^2F(T, X, Y) \begin{bmatrix}
        \vect(\Delta T) \\
        \vect(\Delta X) \\
        \vect(\Delta Y) \\
    \end{bmatrix} = \begin{bmatrix}
        \vect(R_t) \\
        \vect(R_x) \\
        \vect(R_y) \\
    \end{bmatrix},
\end{equation*}
we can perform a similar block elimination as~\eqref{eqn:t-barrier-elim} to show that
\begin{gather*}
    M \begin{bmatrix}
        \vect(\Delta X) \\ \vect(\Delta Y)
    \end{bmatrix} = \begin{bmatrix}
        \vect(R_x) + \grad\log\det(Z)^\top\grad_X P_g(X, Y) \\ \vect(R_y) + \grad\log\det(Z)^\top\grad_Y P_g(X, Y)
    \end{bmatrix} \\
    \Delta T = Z R_t Z + \mathsf{D}_X P_g(X, Y)[\Delta X] + \mathsf{D}_Y P_g(X, Y)[\Delta Y].
\end{gather*}
That is, the main work is in solving linear systems with the matrix $M$, which we do by directly constructing and Cholesky factoring the matrix.

\subsection{Exploiting sparsity}\label{subsec:sparsity}

In many conic problems that arise in practice, the linear constraint matrices $A$ and $G$ are sparse. Where possible, we take advantage of this structure to more efficiently construct the Schur complement matrices $G^\top H G$ or $AH^{-1}A^\top$ required to solve the normal equations~\eqref{eqn:newton-system-elim-b} or~\eqref{eqn:newton-simplified}, respectively, where we recall that $H$ is the Hessian of the barrier function $\grad^2 F(x)$ at some point $x\in\interior\mathcal{K}$. In the following discussion, we focus on constructing $AH^{-1}A^\top$ efficiently where $\mathcal{K}$ consists of a single cone. We remark that the same discussion can be applied to constructing $G^\top H G$, and when $\mathcal{K}$ is a Cartesian product of cones.

For a general cone $\mathcal{K}$, we first construct the matrix $H^{-1}A^\top$ by applying the inverse Hessian-vector product to each row of $A$, then compute $AH^{-1}A^\top$ by performing a matrix multiplication between a sparse and dense matrix. In general, it is difficult to take advantage of sparsity of $A$ when constructing $H^{-1}A^\top$. However, in particular cases, additional structure can be exploited. If $\mathcal{K}$ is the nonnegative orthant, then $H^{-1}$ is a diagonal matrix, and $AH^{-1}A^\top$ can be computed by performing two sparse matrix multiplications. If $\mathcal{K}$ is the positive semidefinite cone, then we follow the method outlined in~\cite{fujisawa1997exploiting}, i.e., we compute each $i,j$-th element of $AH^{-1} A^\top$ by using
\begin{align}\label{eqn:psd-sparse}
    [A H^{-1} A^\top]_{ij} = \sum_{(a,b)\in\mathcal{I}} (A_i)_{ab}^* \biggl( \sum_{(c,d)\in\mathcal{J}}^n (A_j)_{cd} X_{ac} X_{db} \biggr),
\end{align}
where $\mathcal{I}$ and $\mathcal{J}$ are the sets of all nonzero indices in $A_i$ and $A_j$, respectively, i.e., we perform the sum~\eqref{eqn:psd-sparse} only over the nonzero entries of $A_i$ and $A_j$. Note that we only have to compute the upper-triangular elements of $AH^{-1} A^\top$ due to symmetry, and that sorting the rows of $A$ from least to most sparse will reduce the number of floating point operations required. When $A$ has a combination of sparse and dense rows, then we use~\eqref{eqn:psd-sparse} to compute the $i,j$-th element of $AH^{-1} A^\top$ when $A_i$ and $A_j$ are both sparse matrices, and compute all other elements of $A H^{-1} A^\top$ by directly performing the matrix multiplication and inner product $\tr[A_i(X A_j X)]$.

\section{Numerical experiments}\label{sec:experiments}

We compare the performance of QICS to existing state-of-the-art solvers on solving various types of cone programs, with an emphasis on applications arising in quantum information theory. We collect these problems in three new benchmark libraries which we describe in the following sections, and make available at
\begin{center}
    \url{https://github.com/kerry-he/qics-benchmarks}.
\end{center}
Unless otherwise stated, all solvers are run using their default settings. All experiments are run on an Intel i5-11700 CPU with 32GB of RAM. We report the times required to successfully solve problems to \emph{high accuracy}, i.e., relative optimality gap and primal-dual feasibility tolerance of $10^{-8}$, as well as to \emph{low accuracy}, i.e., a relative optimality gap and primal-dual feasibility tolerance of $10^{-5}$. A solver is said to have \emph{failed} if it is unable to reach the desired tolerances within $\SI{3600}{\second}$. In addition to solve times, we also report the normalized shifted geometric mean~\cite{mittelmann}
\begin{equation*}
    \hat{g}_s = \frac{g_s}{\min_s g_s},
\end{equation*}
for solver $s$ to solve a set of problems $\mathcal{P}$, where $g_s$ is the shifted geometric mean defined by
\begin{equation*}
    g_s = \biggl( \prod_{p\in\mathcal{P}} (t_{p,s} + k) \biggr)^{1/\abs{\mathcal{P}}} - k,
\end{equation*}
where $t_{p,s}$ is the time taken to solve problem $p$ using solver $s$, $\abs{\mathcal{P}}$ denotes the number of elements in the set $\mathcal{P}$, and we use a shift of $k=\SI{1}{\second}$. We also present absolute performance profiles as described by~\cite{goulart2024clarabel}, which are defined as the functions
\begin{equation*}
    \rho_s^a(\tau) = \frac{1}{\abs{\mathcal{P}}}  \abs{\{ p \in \mathcal{P} : t_{p,s} \leq \tau \}},
\end{equation*}
and represents the ratio of problems solved by solver $s$ within $\tau$ seconds. Additionally, we present relative performance profiles as described by~\cite{dolan2002benchmarking}, which are defined as the functions
\begin{equation*}
    \rho_s^r(\tau) = \frac{1}{\abs{\mathcal{P}}}  \abs*{ \biggl\{ p \in \mathcal{P} : \frac{t_{p,s}}{\min_s t_{p,s}} \leq \tau \biggr\}},
\end{equation*}
and represents the ratio of problems solved by solver $s$ within a factor of $\tau$ of the fastest time in which the problem was solved.

\subsection{Semidefinite programming}

We first benchmark the performance of QICS on solving sparse semidefinite programs. We compare against state-of-the-art semidefinite programming software, including MOSEK~\cite{mosek}, SDPA~\cite{yamashita2010high,yamashita2012latest}, SDPT3~\cite{toh1999sdpt3}, and SeDuMi~\cite{sturm1999using}. We also benchmark against open-source conic solvers CVXOPT~\cite{vandenberghe2010cvxopt} and Hypatia~\cite{coey2023performance}, upon which our algorithms are heavily based.

\subsubsection{SDPLIB benchmark library}

We first compare these solvers on the SDPLIB library of sparse semidefinite programs~\cite{borchers1999sdplib}. These results are summarized in Table~\ref{tab:sdplib} and Figure~\ref{fig:sdplib}. Overall, these results show that QICS is competitive with state-of-the-art semidefinite programming solvers. Additionally, QICS is much faster than CVXOPT and Hypatia despite sharing similar interior-point algorithms, which we attribute to the fact that CVXOPT and Hypatia\footnote{We use Hypatia's default positive semidefinite cone \texttt{PosSemidefTri}, which does not exploit sparsity. Hypatia provides another cone \texttt{PosSemidefTriSparse} which does exploit sparsity, however users are required to provide the fixed sparsity pattern of the cone a priori. Additionally, the aggregate sparsity patterns of the data matrices for many of the semidefinite programs we benchmark against are not very sparse.} do not exploit sparsity, whereas QICS does. Unsurprisingly, MOSEK is the fastest solver as it is the only commercial, closed-source solver in our benchmarks, and is also fully implemented in C++. Nevertheless, on large graph problems such as \texttt{maxG60}, QICS and MOSEK solve in very similar amounts of time ($\SI{2596}{\second}$ and $\SI{2627}{\second}$, respectively).

\afterpage{%
\begin{table}[t!]
\small
\centering
\caption{Comparison of shifted geometric means and success rates of various semidefinite programming solvers to solve $92$ problems from the SDPLIB benchmark library to full accuracy $\varepsilon=10^{-8}$ and to low accuracy $\varepsilon=10^{-5}$.}
\label{tab:sdplib}
\begin{tabular}{@{}ldddd@{}}
\toprule
 & \multicolumn{2}{c}{\textbf{Full accuracy}} & \multicolumn{2}{c}{\textbf{Low accuracy}} \\ \cmidrule(l{2pt}r{2pt}){2-3}\cmidrule(l{2pt}r{2pt}){4-5} 
 & \mc{Shifted GM} & \mc{Success \%} & \mc{Shifted GM} & \mc{Success \%} \\ \midrule
\textbf{QICS} & 1.76 & 67.4 & 1.81 & 96.7 \\
\textbf{MOSEK} & \multicolumn{1}{Z{.}{.}{-1}}{1.00} & \multicolumn{1}{Z{.}{.}{-1}}{71.7} & \multicolumn{1}{Z{.}{.}{-1}}{1.00} & \multicolumn{1}{Z{.}{.}{-1}@{}}{98.9} \\
\textbf{SDPA} & 10.88 & 42.4 & 7.41 & 81.5 \\
\textbf{SDPT3} & 1.40 & 67.4 & 2.49 & 92.4 \\
\textbf{SeDuMi} & 6.41 & 51.1 & 4.01 & 91.3 \\
\textbf{Hypatia} & 10.96 & 48.9 & 9.19 & 89.1 \\
\textbf{CVXOPT} & 14.71 & 48.9 & 17.11 & 88.0 \\ \bottomrule
\end{tabular}
\end{table}

\begin{figure*}[h!]
\centering
\input{figures/sdplib}
\caption{Comparison of relative performance ratios and solution time profiles of various semidefinite programming solvers to solve $92$ problems from the SDPLIB benchmark library to full accuracy $\varepsilon=10^{-8}$ and to low accuracy $\varepsilon=10^{-5}$.}
\label{fig:sdplib}
\end{figure*}
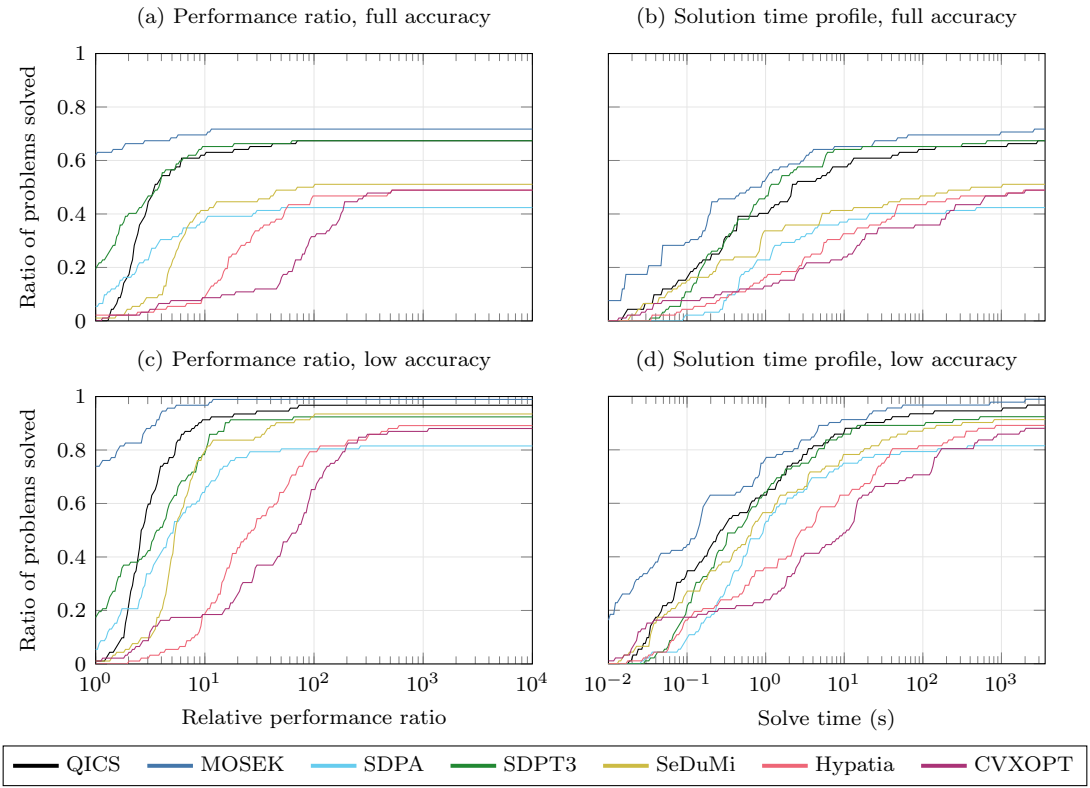
}

\subsubsection{SDPs from quantum information theory}

Next, we compare these solvers on a suite of $49$ semidefinite programs arising from quantum information theory. We briefly describe these applications below.

\begin{itemize}[leftmargin=0em]
    \item[] \textbf{Quantum state discrimination} We design a set of POVMs which maximize the probability of successfully identifying a given set of quantum states. We consider both the minimum-error~\cite[Section 4.2.1]{skrzypczyk2023semidefinite} and unambiguous~\cite[Section 4.2.4]{skrzypczyk2023semidefinite} quantum state discrimination problems.
    
    \item[] \textbf{Quantum state fidelity} The quantum state fidelity $F(X, Y)=\norm{\sqrt{X}\sqrt{Y}}_1^2$ is used to measure how similar two quantum states are. We compute this using the semidefinite program described in~\cite[Section 2.1]{watrous2012simpler}.
    
    \item[] \textbf{Diamond norm} The diamond norm is used to measure the dissimilarity between two quantum channels. This can be computed as the semidefinite program described in~\cite[Section 3.2]{watrous2012simpler}.
    
    \item[] \textbf{Quantum optimal transport} We consider the quantum optimal transport problem~\cite{cole2023quantum}, which minimizes a a linear objective function of a bipartite density matrix subject to constraints on its reduced states. This is analogous to the the classical optimal transport problem, which minimizes a linear objective function of a joint probability distribution subject to satisfying given marginal distributions.
    
    \item[] \textbf{DPS hierarchy} Determining whether a given quantum state is entangled or separable is NP-hard in general. The Doherty-Parrilo-Spedalieri hierarchy~\cite{doherty2004complete} provides a hierarchy of linear matrix inequalities which must be satisfied by all separable states, which we can use to check if a state is separable or entangled (see also, \cite[Section 7]{siddhu2022five}). 
    
    \item[] \textbf{NPA hierarchy} The Navascu\'es-Pironio-Ac\'in hierarchy~\cite{navascues2007bounding,navascues2008convergent} solves non-commutative polynomial optimization problems by constructing a hierarchy of semidefinite programs which lower bound the optimal value. We obtain problems using Moment~\cite{garner2024introducing}, and consider computing bounds for the CHSH inequality, bounds for the I3322 inequality, and the Brown-Fawzi-Fawzi calculation of the conditional quantum entropy~\cite{brown2021device}.
    
    \item[] \textbf{Quantum relative entropy} In~\cite{fawzi2019semidefinite}, it was shown how the quantum relative entropy could be approximated using linear matrix inequalities. We use this technique to solve a toy quantum relative entropy program.

    \item[] \textbf{Quantum chemistry} In~\cite{nakata2001variational,nakata2002density}, it was shown how the potential energy surfaces of molecules could be calculated by solving a semidefinite program. These problems were used to benchmark the performance of SDPA in~\cite{yamashita2003implementation}. We take the specific examples from the SDP benchmarks in~\cite{mittelmann}.

\end{itemize}

\begin{figure*}[t!]
\centering
\input{figures/qsdp}
\caption{Comparison of relative performance ratios and solution time profiles of various semidefinite programming solvers to solve $49$ semidefinite programs arising in quantum information theory to full accuracy $\varepsilon=10^{-8}$ and to low accuracy $\varepsilon=10^{-5}$.}
\label{fig:qsdp}
\end{figure*}
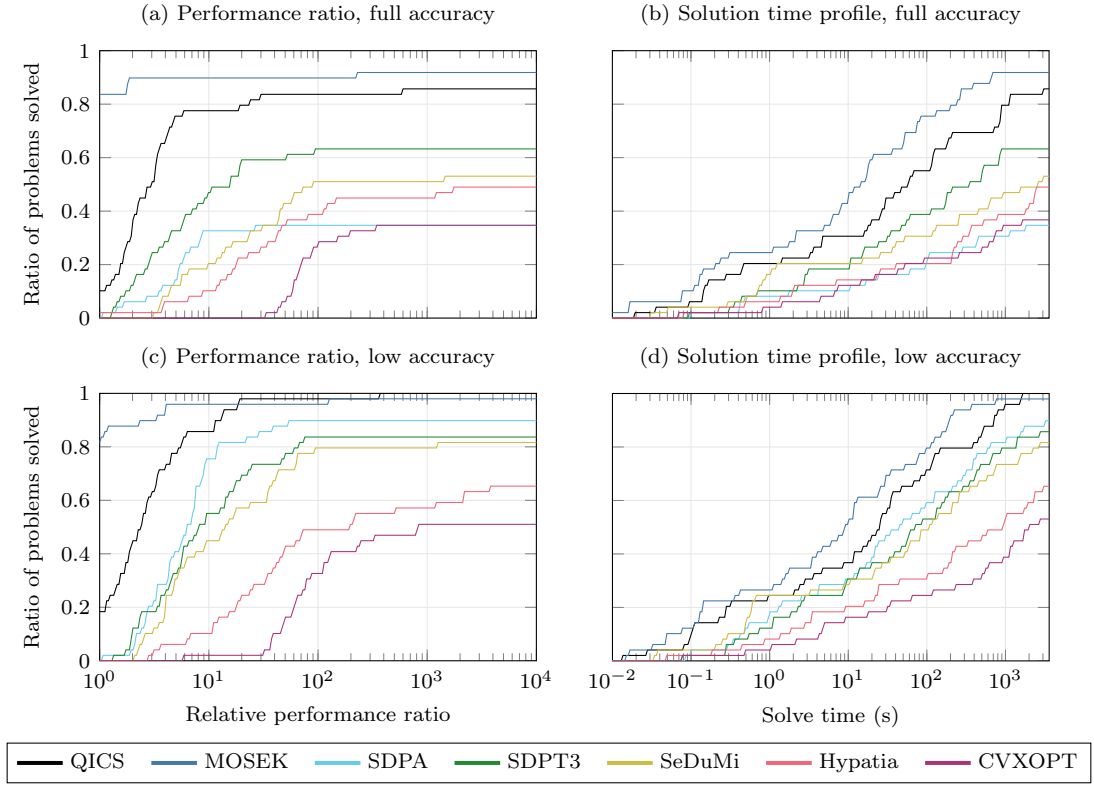

Notably, many of these problems involve semidefinite programs defined on complex Hermitian matrices. To express these problems in a way that solvers which only support real symmetric matrices will accept (i.e., MOSEK, SDPA, SDPT3, and CVXOPT), we use the lifting technique described in~\cite[Section 6.2.7]{aps2020mosek}. The computational results are summarized in Table~\ref{tab:qsdp} and Figure~\ref{fig:qsdp}. Overall, we see similar trends to those observed for the SDPLIB benchmark library. QICS outperforms all other semidefinite programming solvers except MOSEK, which it is also relatively competitive with. Additionally, QICS is able to solve all problems to low accuracy, which all other solvers are unable to do. As this suite of problems contained more large-scale problems compared to the SDPLIB benchmark library, timeout and insufficient memory errors were more common, particularly for the CVXOPT and Hypatia solvers.

\begin{table}[t!]
\small
\centering
\caption{Comparison of shifted geometric means and success rates of various semidefinite programming solvers to solve $49$ semidefinite programs arising in quantum information theory to full accuracy $\varepsilon=10^{-8}$ and to low accuracy $\varepsilon=10^{-5}$.}
\label{tab:qsdp}
\begin{tabular}{@{}ldddd@{}}
\toprule
 & \multicolumn{2}{c}{\textbf{Full accuracy}} & \multicolumn{2}{c}{\textbf{Low accuracy}} \\ \cmidrule(l{2pt}r{2pt}){2-3}\cmidrule(l{2pt}r{2pt}){4-5} 
 & \mc{Shifted GM} & \mc{Success \%} & \mc{Shifted GM} & \mc{Success \%} \\ \midrule
\textbf{QICS} & 2.75 & 85.7 & 1.89 & \multicolumn{1}{Z{.}{.}{-1}@{}}{100.0} \\
\textbf{MOSEK} & \multicolumn{1}{Z{.}{.}{-1}}{1.00} & \multicolumn{1}{Z{.}{.}{-1}}{91.8} & \multicolumn{1}{Z{.}{.}{-1}}{1.00} & 98.0 \\
\textbf{SDPA} & 29.29 & 34.7 & 4.08 & 89.8 \\
\textbf{SDPT3} & 8.41 & 63.3 & 5.59 & 85.7 \\
\textbf{SeDuMi} & 12.97 & 53.1 & 6.26 & 81.6 \\
\textbf{Hypatia} & 23.30 & 49.0 & 18.88 & 65.3 \\
\textbf{CVXOPT} & 30.35 & 36.7 & 34.57 & 53.1 \\ \bottomrule
\end{tabular}
\end{table}

\subsection{Quantum relative entropy programs}\label{subsec:exp-qrep}

We also benchmark QICS on solving quantum relative entropy programs. We compare against Hypatia~\cite{coey2023performance}, DDS~\cite{karimi2023efficient,karimi2024domain}, as well as modelling the problems with CVXQUAD~\cite{fawzi2019semidefinite} and solving the resulting semideifnite program using MOSEK. We also compare between solving problems using the full suite of cones supported by QICS, including the quantum conditional entropy and quantum key distribution cones, against solving problems only using the quantum relative entropy cone. We denote these latter results with an asterisk, i.e., QICS*. To compare these solvers, we collected a suite of $144$ quantum relative entropy programs from across a range of problems arising from quantum information theory. We briefly describe these applications below.

\begin{itemize}[leftmargin=0em]

    \item[] \textbf{Nearest correlation matrix}
    We compute the closest correlation matrix (i.e., symmetric, positive semidefinite, and diagonal entries are all one) to a given matrix, in the quantum relative entropy sense (see, e.g.,~\cite[Section 8.1]{karimi2023efficient}). We also consider variations of these problems where the matrix is restricted to being tridiagonal, as considered in~\cite{karimi2023efficient}.
    
    \item[] \textbf{Relative entropy of entanglement}
    We quantify how entangled a given bipartite quantum state is by estimating the distance between the matrix and the set of matrices satisfying the positive partial transpose criteria~\cite{horodecki1996separability} (see, e.g.,~\cite[Section 4]{fawzi2018efficient}).
    
    \item[] \textbf{Quantum channel capacities}
    Quantum channel capacities quantify the maximum amount of information we can transit through a noisy quantum channel. We consider the classical-quantum channel capacity (see, e.g.,~\cite[Section 8.4.4]{coey2023conic}), the entanglement-assisted channel capacity (see, e.g.,~\cite[Section 3.3.1]{he2024exploiting}), and the quantum-quantum channel capacity of degradable channels (see, e.g.,~\cite[Section 3.3.2]{he2024exploiting}). Note that the classical-quantum capacity can be modeled using just the quantum entropy cone, and the entanglement-assisted and quantum-quantum capacities can be modeled using just the quantum conditional entropy cone.
    
    \item[] \textbf{Quantum rate distortion}
    The quantum rate distortion function characterizes the minimum amount of information required to compress and transmit a quantum information source through a quantum channel without exceeding a given distortion threshold (see, e.g.,~\cite[Section 3.2]{he2024exploiting}). Note that we can model this problem using the quantum conditional entropy cone. We also study the entanglement-assisted rate distortion function, which we can simplify using symmetry reduction as described in~\cite{he2024efficient}, and which we can model using a Cartesian product of the classical relative entropy cone and quantum relative entropy cone, as described in~\cite[Section 6.2.1]{he2024exploiting}.
    
    \item[] \textbf{Quantum key rates}
    The quantum key rate is a quantity which characterizes the security of a given quantum protocol (see, e.g.,~\cite[Section 2]{winick2018reliable}). This problem can be modelled using the quantum key distribution cone. We consider specific problem quantum key distribution protocols from~\cite{hu2022robust,lorente2024quantum}.
    
    \item[] \textbf{Ground energy of Hamiltonian}
    We lower bound the ground energy of a Hamiltonian using the method in~\cite{fawzi2023entropy} (see also, \cite[Section 3.4]{he2024exploiting}). Note that this problem can modelled using the quantum conditional entropy cone.
\end{itemize}

\begin{figure*}[t!]
\centering
\input{figures/qreps}
\caption{Comparison of relative performance ratios and solution time profiles of various solvers to solve $144$ quantum relative entropy programs to full accuracy $\varepsilon=10^{-8}$ and to low accuracy $\varepsilon=10^{-5}$. Note that QICS refers to results using the full suite of cones we implement, whereas QICS* refers to results using only the quantum entropy and quantum relative entropy cones.}
\label{fig:qreps}
\end{figure*}
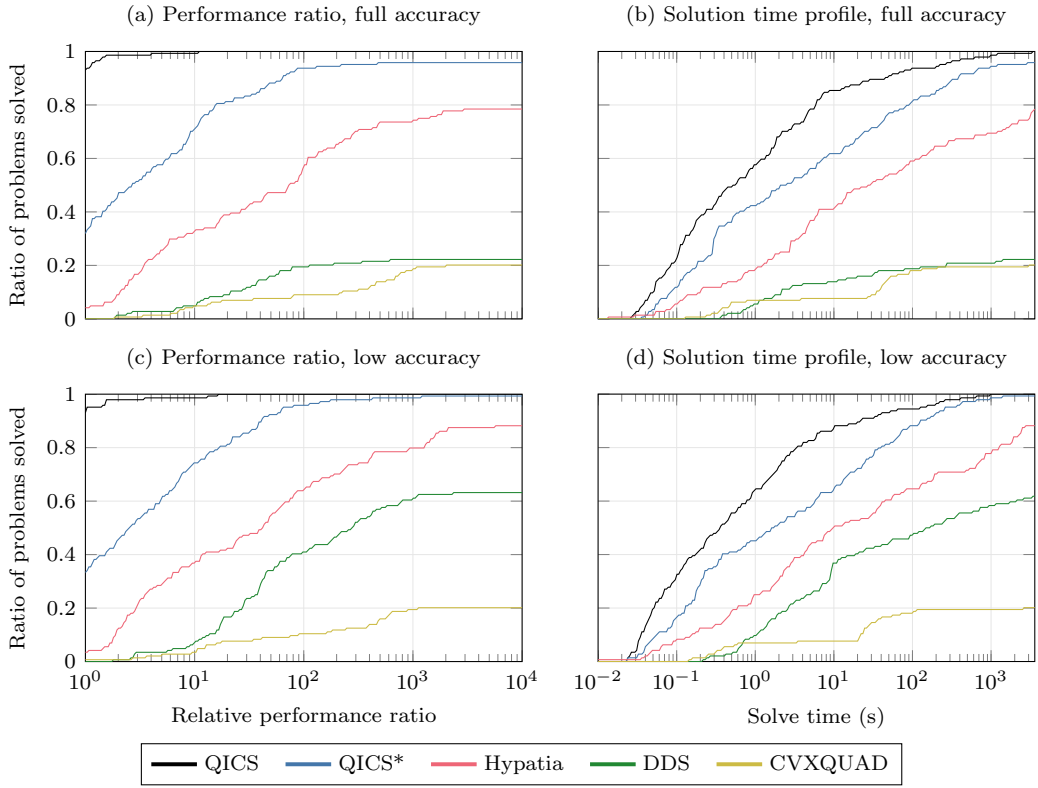

These results are summarized in Table~\ref{tab:qreps} and Figure~\ref{fig:qreps}. Overall, the results show that QICS is faster than all existing quantum relative entropy programming software currently available. Additionally, using the full suite of cones suuported by QICS provides a further order of magnitude improvement in solve times over just using the quantum relative entropy cone, while also improving the solver's success rate to 100\%. Despite our stepping algorithm being very similar to Hypatia, our solver has more efficient implementations of certain key steps of the algorithm, such as having more efficient cone oracles. DDS has relatively competitive solve times per iteration, however we found that DDS often required more iterations to converge, and was unable to converge to high accuracies for many problems. As expected, CVXQUAD was also only able to successfully solve problems of relatively small dimensions due to using linear matrix inequalities with $2n^2\times2n^2$ matrices to approximate quantum relative entropy involving $n\times n$ matrix arguments.

\begin{table}[t!]
\small
\centering
\caption{Comparison of shifted geometric means and success rates of various solvers to solve $144$ quantum relative entropy programs to full accuracy $\varepsilon=10^{-8}$ and to low accuracy $\varepsilon=10^{-5}$. Note that QICS refers to results using the full suite of cones we implement, whereas QICS* refers to results using on the quantum entropy and quantum relative entropy cones.}
\label{tab:qreps}
\begin{tabular}{@{}ldddd@{}}
\toprule
 & \multicolumn{2}{c}{\textbf{Full accuracy}} & \multicolumn{2}{c}{\textbf{Low accuracy}} \\ \cmidrule(l{2pt}r{2pt}){2-3}\cmidrule(l{2pt}r{2pt}){4-5} 
 & \mc{Shifted GM} & \mc{Success \%} & \mc{Shifted GM} & \mc{Success \%} \\ \midrule
\textbf{QICS} & \multicolumn{1}{Z{.}{.}{-1}}{1.00} & \multicolumn{1}{Z{.}{.}{-1}}{100.0} & \multicolumn{1}{Z{.}{.}{-1}}{1.00} & \multicolumn{1}{Z{.}{.}{-1}@{}}{100.0} \\
\textbf{QICS*} & 3.68 & 95.8 & 3.24 & 99.3 \\
\textbf{Hypatia} & 23.34 & 78.5 & 18.64 & 88.2 \\
\textbf{DDS} & 428.12 & 22.2 & 75.08 & 61.2 \\
\textbf{CVXQUAD} & 539.47 & 20.1 & 704.78 & 20.1 \\ \bottomrule
\end{tabular}
\end{table}

\subsection{Noncommutative perspective}

Finally, we compare the performance of QICS and CVXQUAD using MOSEK when solving convex optimization problems involving noncommutative perspectives. We include results for both QICS using the technique described in Section~\ref{eqn:avoid-inv-hess} to avoid inverse Hessian-vector products where appropriate, and QICS without using this technique, which we denote with an asterisk, i.e., QICS*. To compare these solvers, we benchmark against a total of $48$ problems drawn from the following applications.

\begin{itemize}[leftmargin=0em]
    \item[] \textbf{D-optimal design} In~\cite[Section 4.4]{bertsimas2023new}, a matrix perspective reformulation technique is proposed to obtain a relaxation of the D-optimal design problem by minimizing the trace of an operator relative entropy. This was shown to achieve empirically tighter bounds compared to other relaxation strategies.
    
    \item[] \textbf{Measured relative entropy} Measured relative entropies are functions used to measure the amount of dissimilarity between two quantum states which arise in quantum hypothesis testing tasks. Variational expressions were given for these functions in~\cite{berta2017variational}, and it was shown in~\cite[Propositions 3 and 7]{huang2024semi} how these expressions could be given as the solution to conic programs which optimize over the epigraph of noncommutative perspective functions. We focus on computing the measured R\'enyi relative entropy of states, which can be computed by optimizing over the weighted matrix geometric mean.
    
    \item[] \textbf{Expectation values of equilibrium states} In quantum many-body theory, an important problem is to describe observables in equilibrium states. In~\cite{fawzi2024certified}, it was shown how a hierarchy of conic problems involving the epigraph of the operator relative entropy can produce converging bounds to the expectation values of these observables. 
\end{itemize}

\begin{table}[t!]
\small
\centering
\caption{Comparison of shifted geometric means and success rates of various solvers to solve $48$ conic programs involving noncommutative perspective functions to full accuracy $\varepsilon=10^{-8}$ and to low accuracy $\varepsilon=10^{-5}$. Note that QICS refers to results employing the technique described in Section~\ref{eqn:avoid-inv-hess} where appropriate, and QICS* refers to results without using this technique.}
\label{tab:ncps}
\begin{tabular}{@{}ldddd@{}}
\toprule
 & \multicolumn{2}{c}{\textbf{Full accuracy}} & \multicolumn{2}{c}{\textbf{Low accuracy}} \\ \cmidrule(l{2pt}r{2pt}){2-3}\cmidrule(l{2pt}r{2pt}){4-5} 
 & \mc{Shifted GM} & \mc{Success \%} & \mc{Shifted GM} & \mc{Success \%} \\ \midrule
\textbf{QICS} & \multicolumn{1}{Z{.}{.}{-1}}{1.00} & 83.3 & \multicolumn{1}{Z{.}{.}{-1}}{1.00} & \multicolumn{1}{Z{.}{.}{-1}@{}}{100.0} \\
\textbf{QICS*} & 1.69 & 79.2 & 1.42 & 87.5 \\
\textbf{CVXQUAD} & 1.06 & \multicolumn{1}{Z{.}{.}{-1}}{87.5} & 1.44 & 87.5 \\ \bottomrule
\end{tabular}
\end{table}

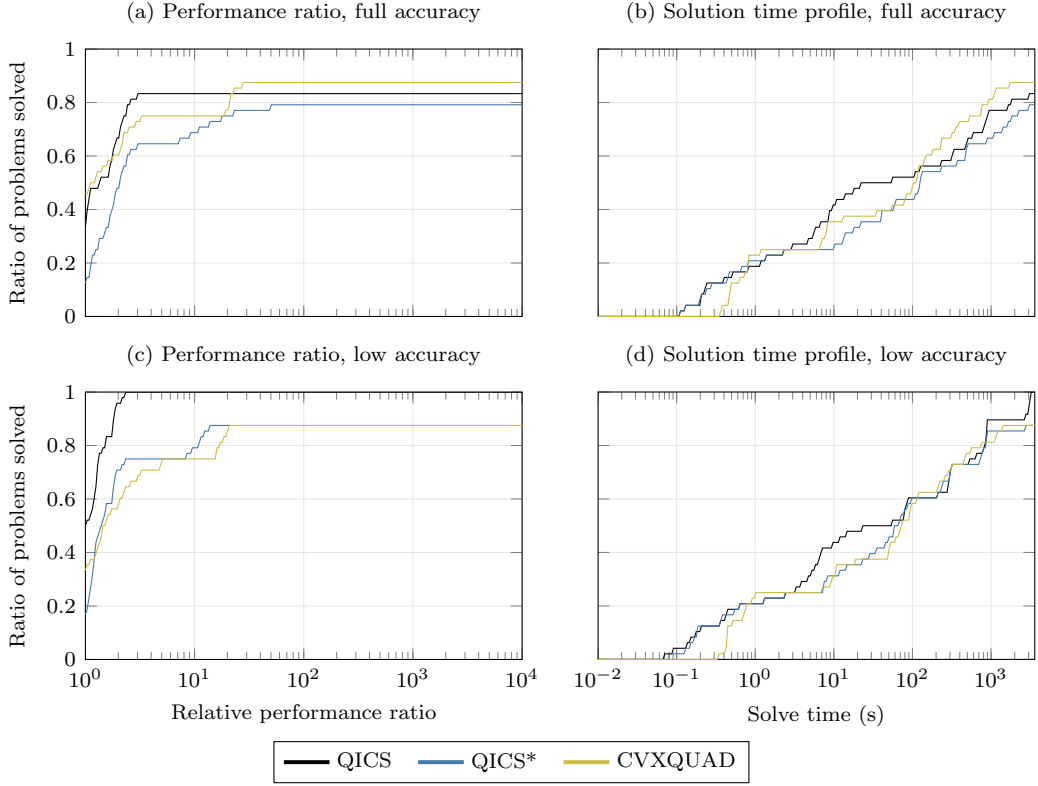
\begin{figure*}[t!]
\centering
\input{figures/ncps}
\caption{Comparison of relative performance ratios and solution time profiles of various solvers to solve $48$ conic programs involving noncommutative perspective functions to full accuracy $\varepsilon=10^{-8}$ and to low accuracy $\varepsilon=10^{-5}$. Note that QICS refers to results employing the technique described in Section~\ref{eqn:avoid-inv-hess} where appropriate, and QICS* refers to results without using this technique.}
\label{fig:ncps}
\end{figure*}

These results are summarized in Table~\ref{tab:ncps} and Figure~\ref{fig:ncps}. Overall, we observe that QICS* and CVXQUAD have fairly comparable performance, where CVXQUAD is generally faster and more stable when solving problems to high accuracy, and QICS* and CVXQUAD perform almost identically when solving to low accuracy. We also see that avoiding inverse Hessian-vector products using the technique from Section~\ref{eqn:avoid-inv-hess} results in a notable improvement in computation times. This is particularly the case for computing the expectation values of equilibrium states, which minimizes a $2^L\times2^L$ decision variable over the epigraph of the noncommutative perspective involving significantly larger $2^{L-2}\times2^{L-2}$ matrices. By avoiding solving linear systems with the Hessian of the barrier function of this epigraph, we reduce the work from building and Cholesky factoring a $2^{4L-3}\times2^{4L-3}$ Hessian matrix to constructing and Cholesky factoring a $2^{2L}\times2^{2L}$ matrix. This allows QICS to solve problems of this type for up to $L=5$, whereas QICS* and CVXQUAD were unable to solve these problems at $L=5$ due to insufficient memory.




\section{Conclusion}

We have presented QICS, an open-source conic solver tailored towards problems arising in quantum information theory. Our benchmarks show that QICS can solve quantum relative entropy programs significantly faster than existing software, and can solve semidefinite programs at comparable speeds to existing state-of-the-art semidefinite programming solvers. We outline some directions for future work below.

\paragraph{Presolver}
Many state-of-the-art linear and semidefinite programming software have a presolver which aims to simplify the given problem before passing it to the interior point algorithm, see, e.g.,~\cite{andersen1995presolving}. It would be interesting to see how to tailor these ideas for quantum relative entropy programming. For example, it is not always obvious what the most efficient way to model a given conic program is, similar to how it may be more efficient to solve the primal or dual of a given semidefinite program. Additionally, it would be helpful to automatically detect when we can reduce a problem from using the full quantum relative entropy cone, to a cone with more efficient oracles such as the quantum conditional entropy cone. Similarly, facial reduction preprocessing steps~\cite{drusvyatskiy2017many} would be useful for many problems arising from quantum information theory.

\paragraph{Large-scale problems}
Currently, QICS is tailored towards solving moderately sizes problems to high accuracy. However, as problems in quantum information theory scale exponentially with the number of qubits involved, it is of practical interest to solve large scale quantum relative entropy programs. One way to do this is to explore first-order methods, such as those discussed in~\cite{he2024bregman,he2024exploiting}. Alternatively, we can augment our interior point algorithm with techniques such as preconditioned conjugate gradient to solve the linear systems that arise in the algorithm, which allows us to avoid having to explicitly form large Hessian or Schur complement matrices (see, e.g.,~\cite{wolkowicz2004solving,tang2024feasible}).

\section*{Acknowledgements}

K.\ He was supported by an Australian Government Research Training Program (RTP) Scholarship. H.\ Fawzi was partially funded by UK Research and Innovation (UKRI) under the UK government’s Horizon Europe funding guarantee EP/X032051/1. This research was also supported, in part, by an Australian Research Council Discovery Project (number DP250104201), funded by the Australian Government.

\bibliographystyle{IEEEtran}
\bibliography{refs}

\end{document}

%% file: figures/sdplib.tex
\pgfplotstableread[col sep=comma]{figures/data/sdplib.csv}\data
\begin{tikzpicture}

    \definecolor{clr1}{RGB}{68,119,170}
    \definecolor{clr2}{RGB}{102,204,238}
    \definecolor{clr3}{RGB}{34,136,51}
    \definecolor{clr4}{RGB}{204,187,68}
    \definecolor{clr5}{RGB}{238,102,119}
    \definecolor{clr6}{RGB}{170,51,119}
    \definecolor{clr7}{RGB}{0,0,0}

    \begin{semilogxaxis}[
        name=mainplot,
        width = 0.5\textwidth,
        height = 0.35\textwidth,
        xmin = 1, xmax = 10000,
        ymin = 0, ymax=1,
        ylabel=Ratio of problems solved,
        label style = {font=\footnotesize},
        ticklabel style = {font=\footnotesize},
        title style = {font=\footnotesize},
        no markers,
        grid=major,
        major grid style={line width=.1pt,draw=gray!20},
        yminorticks=false,
        xticklabels={,,},
        title = {(a) Performance ratio, full accuracy}
        ]
        
        \addplot[clr7] table[x=tau-r , y=ours-rf] {\data}; 
        \addplot[clr1] table[x=tau-r , y=mosek-rf] {\data};
        \addplot[clr2] table[x=tau-r , y=sdpa-rf] {\data};
        \addplot[clr3] table[x=tau-r , y=sdpt3-rf] {\data};
        \addplot[clr4] table[x=tau-r , y=sedumi-rf] {\data};
        \addplot[clr5] table[x=tau-r , y=hypatia-rf] {\data};
        \addplot[clr6] table[x=tau-r , y=cvxopt-rf] {\data};

    \end{semilogxaxis}

    \begin{semilogxaxis}[
        name=secondplot,
        at={(mainplot.north east)},
        xshift=1cm,
        anchor=north west,   
        width = 0.5\textwidth,
        height = 0.35\textwidth,
        xmin = 1e-2, xmax = 3600,
        ymin = 0, ymax=1,
        label style = {font=\footnotesize},
        ticklabel style = {font=\footnotesize},
        title style = {font=\footnotesize},
        no markers,
        grid=major,
        major grid style={line width=.1pt,draw=gray!20},
        yminorticks=false,
        yticklabels={,,},
        xticklabels={,,},
        title = {(b) Solution time profile, full accuracy}
        ]
        
        \addplot[clr7] table[x=tau-t , y=ours-tf] {\data}; 
        \addplot[clr1] table[x=tau-t , y=mosek-tf] {\data};
        \addplot[clr2] table[x=tau-t , y=sdpa-tf] {\data};
        \addplot[clr3] table[x=tau-t , y=sdpt3-tf] {\data};
        \addplot[clr4] table[x=tau-t , y=sedumi-tf] {\data};
        \addplot[clr5] table[x=tau-t , y=hypatia-tf] {\data};
        \addplot[clr6] table[x=tau-t , y=cvxopt-tf] {\data};
        
    \end{semilogxaxis}

    \begin{semilogxaxis}[
        name=thirdplot,
        at={(mainplot.south west)},
        yshift=-1cm,
        anchor=north west,   
        width = 0.5\textwidth,
        height = 0.35\textwidth,
        xmin = 1, xmax = 10000,
        ymin = 0, ymax=1,
        xlabel=Relative performance ratio,
        ylabel=Ratio of problems solved,
        label style = {font=\footnotesize},
        ticklabel style = {font=\footnotesize},
        title style = {font=\footnotesize},
        no markers,
        grid=major,
        major grid style={line width=.1pt,draw=gray!20},
        yminorticks=false,
        title = {(c) Performance ratio, low accuracy}
        ]
        
        \addplot[clr7] table[x=tau-r , y=ours-rn] {\data}; 
        \addplot[clr1] table[x=tau-r , y=mosek-rn] {\data};
        \addplot[clr2] table[x=tau-r , y=sdpa-rn] {\data};
        \addplot[clr3] table[x=tau-r , y=sdpt3-rn] {\data};
        \addplot[clr4] table[x=tau-r , y=sedumi-rn] {\data};
        \addplot[clr5] table[x=tau-r , y=hypatia-rn] {\data};
        \addplot[clr6] table[x=tau-r , y=cvxopt-rn] {\data};
        
    \end{semilogxaxis}

    \begin{semilogxaxis}[
        at={(thirdplot.north east)},
        xshift=1cm,
        anchor=north west,   
        width = 0.5\textwidth,
        height = 0.35\textwidth,
        xmin = 1e-2, xmax = 3600,
        ymin = 0, ymax=1,
        xlabel=Solve time (s),
        label style = {font=\footnotesize},
        ticklabel style = {font=\footnotesize},
        title style = {font=\footnotesize},
        no markers,
        grid=major,
        major grid style={line width=.1pt,draw=gray!20},
        yminorticks=false,
        yticklabels={,,},
        title = {(d) Solution time profile, low accuracy},
        legend cell align={left},
        legend columns=-1,
        legend style={/tikz/every even column/.append style={column sep=0.25cm}, font=\footnotesize, at={(1.1,-0.3)},cells={line width=1.5pt}},
        ]
        
        \addplot[clr7] table[x=tau-t , y=ours-tn] {\data}; 
        \addplot[clr1] table[x=tau-t , y=mosek-tn] {\data};
        \addplot[clr2] table[x=tau-t , y=sdpa-tn] {\data};
        \addplot[clr3] table[x=tau-t , y=sdpt3-tn] {\data};
        \addplot[clr4] table[x=tau-t , y=sedumi-tn] {\data};
        \addplot[clr5] table[x=tau-t , y=hypatia-tn] {\data};
        \addplot[clr6] table[x=tau-t , y=cvxopt-tn] {\data};

        \legend{
            QICS,
            MOSEK,
            SDPA,
            SDPT3,
            SeDuMi,
            Hypatia,
            CVXOPT
        }
        
    \end{semilogxaxis}

\end{tikzpicture}

%% file: figures/qsdp.tex
\pgfplotstableread[col sep=comma]{figures/data/qsdp.csv}\data
\begin{tikzpicture}

    \definecolor{clr1}{RGB}{68,119,170}
    \definecolor{clr2}{RGB}{102,204,238}
    \definecolor{clr3}{RGB}{34,136,51}
    \definecolor{clr4}{RGB}{204,187,68}
    \definecolor{clr5}{RGB}{238,102,119}
    \definecolor{clr6}{RGB}{170,51,119}
    \definecolor{clr7}{RGB}{0,0,0}

    \begin{semilogxaxis}[
        name=mainplot,
        width = 0.5\textwidth,
        height = 0.35\textwidth,
        xmin = 1, xmax = 10000,
        ymin = 0, ymax=1,
        ylabel=Ratio of problems solved,
        label style = {font=\footnotesize},
        ticklabel style = {font=\footnotesize},
        title style = {font=\footnotesize},
        no markers,
        grid=major,
        major grid style={line width=.1pt,draw=gray!20},
        yminorticks=false,
        xticklabels={,,},
        title = {(a) Performance ratio, full accuracy}
        ]
        
        \addplot[clr7] table[x=tau-r , y=ours-rf] {\data}; 
        \addplot[clr1] table[x=tau-r , y=mosek-rf] {\data};
        \addplot[clr2] table[x=tau-r , y=sdpa-rf] {\data};
        \addplot[clr3] table[x=tau-r , y=sdpt3-rf] {\data};
        \addplot[clr4] table[x=tau-r , y=sedumi-rf] {\data};
        \addplot[clr5] table[x=tau-r , y=hypatia-rf] {\data};
        \addplot[clr6] table[x=tau-r , y=cvxopt-rf] {\data};

    \end{semilogxaxis}

    \begin{semilogxaxis}[
        name=secondplot,
        at={(mainplot.north east)},
        xshift=1cm,
        anchor=north west,   
        width = 0.5\textwidth,
        height = 0.35\textwidth,
        xmin = 1e-2, xmax = 3600,
        ymin = 0, ymax=1,
        label style = {font=\footnotesize},
        ticklabel style = {font=\footnotesize},
        title style = {font=\footnotesize},
        no markers,
        grid=major,
        major grid style={line width=.1pt,draw=gray!20},
        yminorticks=false,
        yticklabels={,,},
        xticklabels={,,},
        title = {(b) Solution time profile, full accuracy}
        ]
        
        \addplot[clr7] table[x=tau-t , y=ours-tf] {\data}; 
        \addplot[clr1] table[x=tau-t , y=mosek-tf] {\data};
        \addplot[clr2] table[x=tau-t , y=sdpa-tf] {\data};
        \addplot[clr3] table[x=tau-t , y=sdpt3-tf] {\data};
        \addplot[clr4] table[x=tau-t , y=sedumi-tf] {\data};
        \addplot[clr5] table[x=tau-t , y=hypatia-tf] {\data};
        \addplot[clr6] table[x=tau-t , y=cvxopt-tf] {\data};
        
    \end{semilogxaxis}

    \begin{semilogxaxis}[
        name=thirdplot,
        at={(mainplot.south west)},
        yshift=-1cm,
        anchor=north west,   
        width = 0.5\textwidth,
        height = 0.35\textwidth,
        xmin = 1, xmax = 10000,
        ymin = 0, ymax=1,
        xlabel=Relative performance ratio,
        ylabel=Ratio of problems solved,
        label style = {font=\footnotesize},
        ticklabel style = {font=\footnotesize},
        title style = {font=\footnotesize},
        no markers,
        grid=major,
        major grid style={line width=.1pt,draw=gray!20},
        yminorticks=false,
        title = {(c) Performance ratio, low accuracy}
        ]
        
        \addplot[clr7] table[x=tau-r , y=ours-rn] {\data}; 
        \addplot[clr1] table[x=tau-r , y=mosek-rn] {\data};
        \addplot[clr2] table[x=tau-r , y=sdpa-rn] {\data};
        \addplot[clr3] table[x=tau-r , y=sdpt3-rn] {\data};
        \addplot[clr4] table[x=tau-r , y=sedumi-rn] {\data};
        \addplot[clr5] table[x=tau-r , y=hypatia-rn] {\data};
        \addplot[clr6] table[x=tau-r , y=cvxopt-rn] {\data};
        
    \end{semilogxaxis}

    \begin{semilogxaxis}[
        at={(thirdplot.north east)},
        xshift=1cm,
        anchor=north west,   
        width = 0.5\textwidth,
        height = 0.35\textwidth,
        xmin = 1e-2, xmax = 3600,
        ymin = 0, ymax=1,
        xlabel=Solve time (s),
        label style = {font=\footnotesize},
        ticklabel style = {font=\footnotesize},
        title style = {font=\footnotesize},
        no markers,
        grid=major,
        major grid style={line width=.1pt,draw=gray!20},
        yminorticks=false,
        yticklabels={,,},
        title = {(d) Solution time profile, low accuracy},
        legend cell align={left},
        legend columns=-1,
        legend style={/tikz/every even column/.append style={column sep=0.25cm}, font=\footnotesize, at={(1.1,-0.3)},cells={line width=1.5pt}},
        ]
        
        \addplot[clr7] table[x=tau-t , y=ours-tn] {\data}; 
        \addplot[clr1] table[x=tau-t , y=mosek-tn] {\data};
        \addplot[clr2] table[x=tau-t , y=sdpa-tn] {\data};
        \addplot[clr3] table[x=tau-t , y=sdpt3-tn] {\data};
        \addplot[clr4] table[x=tau-t , y=sedumi-tn] {\data};
        \addplot[clr5] table[x=tau-t , y=hypatia-tn] {\data};        
        \addplot[clr6] table[x=tau-t , y=cvxopt-tn] {\data};

        \legend{
            QICS,
            MOSEK,
            SDPA,
            SDPT3,
            SeDuMi,
            Hypatia,
            CVXOPT
        }
        
    \end{semilogxaxis}

\end{tikzpicture}

%% file: figures/qreps.tex
\pgfplotstableread[col sep=comma]{figures/data/qreps.csv}\data
\begin{tikzpicture}

    \definecolor{clr1}{RGB}{68,119,170}
    \definecolor{clr2}{RGB}{102,204,238}
    \definecolor{clr3}{RGB}{34,136,51}
    \definecolor{clr4}{RGB}{204,187,68}
    \definecolor{clr5}{RGB}{238,102,119}
    \definecolor{clr6}{RGB}{170,51,119}
    \definecolor{clr7}{RGB}{0,0,0}

    \begin{semilogxaxis}[
        name=mainplot,
        width = 0.5\textwidth,
        height = 0.35\textwidth,
        xmin = 1, xmax = 10000,
        ymin = 0, ymax=1,
        ylabel=Ratio of problems solved,
        label style = {font=\footnotesize},
        ticklabel style = {font=\footnotesize},
        title style = {font=\footnotesize},
        no markers,
        grid=major,
        major grid style={line width=.1pt,draw=gray!20},
        yminorticks=false,
        xticklabels={,,},
        title = {(a) Performance ratio, full accuracy}
        ]
        
        \addplot[clr7] table[x=tau-r , y=qics-rf] {\data}; 
        \addplot[clr1] table[x=tau-r , y=qics*-rf] {\data};
        \addplot[clr5] table[x=tau-r , y=hypatia-rf] {\data};
        \addplot[clr3] table[x=tau-r , y=dds-rf] {\data};
        \addplot[clr4] table[x=tau-r , y=cvxquad-rf] {\data};

    \end{semilogxaxis}

    \begin{semilogxaxis}[
        name=secondplot,
        at={(mainplot.north east)},
        xshift=1cm,
        anchor=north west,   
        width = 0.5\textwidth,
        height = 0.35\textwidth,
        xmin = 1e-2, xmax = 3600,
        ymin = 0, ymax=1,
        label style = {font=\footnotesize},
        ticklabel style = {font=\footnotesize},
        title style = {font=\footnotesize},
        no markers,
        grid=major,
        major grid style={line width=.1pt,draw=gray!20},
        yminorticks=false,
        yticklabels={,,},
        xticklabels={,,},
        title = {(b) Solution time profile, full accuracy}
        ]
        
        \addplot[clr7] table[x=tau-t , y=qics-tf] {\data}; 
        \addplot[clr1] table[x=tau-t , y=qics*-tf] {\data};
        \addplot[clr5] table[x=tau-t , y=hypatia-tf] {\data};
        \addplot[clr3] table[x=tau-t , y=dds-tf] {\data};
        \addplot[clr4] table[x=tau-t , y=cvxquad-tf] {\data};
        
    \end{semilogxaxis}

    \begin{semilogxaxis}[
        name=thirdplot,
        at={(mainplot.south west)},
        yshift=-1cm,
        anchor=north west,   
        width = 0.5\textwidth,
        height = 0.35\textwidth,
        xmin = 1, xmax = 10000,
        ymin = 0, ymax=1,
        xlabel=Relative performance ratio,
        ylabel=Ratio of problems solved,
        label style = {font=\footnotesize},
        ticklabel style = {font=\footnotesize},
        title style = {font=\footnotesize},
        no markers,
        grid=major,
        major grid style={line width=.1pt,draw=gray!20},
        yminorticks=false,
        title = {(c) Performance ratio, low accuracy}
        ]
        
        \addplot[clr7] table[x=tau-r , y=qics-rn] {\data}; 
        \addplot[clr1] table[x=tau-r , y=qics*-rn] {\data};
        \addplot[clr5] table[x=tau-r , y=hypatia-rn] {\data};
        \addplot[clr3] table[x=tau-r , y=dds-rn] {\data};
        \addplot[clr4] table[x=tau-r , y=cvxquad-rn] {\data};
        
    \end{semilogxaxis}

    \begin{semilogxaxis}[
        at={(thirdplot.north east)},
        xshift=1cm,
        anchor=north west,   
        width = 0.5\textwidth,
        height = 0.35\textwidth,
        xmin = 1e-2, xmax = 3600,
        ymin = 0, ymax=1,
        xlabel=Solve time (s),
        label style = {font=\footnotesize},
        ticklabel style = {font=\footnotesize},
        title style = {font=\footnotesize},
        no markers,
        grid=major,
        major grid style={line width=.1pt,draw=gray!20},
        yminorticks=false,
        yticklabels={,,},
        title = {(d) Solution time profile, low accuracy},
        legend cell align={left},
        legend columns=-1,
        legend style={/tikz/every even column/.append style={column sep=0.25cm}, font=\footnotesize, at={(0.7,-0.3)},cells={line width=1.5pt}},
        ]
        
        \addplot[clr7] table[x=tau-t , y=qics-tn] {\data}; 
        \addplot[clr1] table[x=tau-t , y=qics*-tn] {\data};
        \addplot[clr5] table[x=tau-t , y=hypatia-tn] {\data};
        \addplot[clr3] table[x=tau-t , y=dds-tn] {\data};
        \addplot[clr4] table[x=tau-t , y=cvxquad-tn] {\data};

        \legend{
            QICS,
            QICS*,
            Hypatia,
            DDS,
            CVXQUAD
        }
        
    \end{semilogxaxis}

\end{tikzpicture}

%% file: figures/ncps.tex
\pgfplotstableread[col sep=comma]{figures/data/ncps.csv}\data
\begin{tikzpicture}

    \definecolor{clr1}{RGB}{68,119,170}
    \definecolor{clr2}{RGB}{102,204,238}
    \definecolor{clr3}{RGB}{34,136,51}
    \definecolor{clr4}{RGB}{204,187,68}
    \definecolor{clr5}{RGB}{238,102,119}
    \definecolor{clr6}{RGB}{170,51,119}
    \definecolor{clr7}{RGB}{0,0,0}

    \begin{semilogxaxis}[
        name=mainplot,
        width = 0.5\textwidth,
        height = 0.35\textwidth,
        xmin = 1, xmax = 10000,
        ymin = 0, ymax=1,
        ylabel=Ratio of problems solved,
        label style = {font=\footnotesize},
        ticklabel style = {font=\footnotesize},
        title style = {font=\footnotesize},
        no markers,
        grid=major,
        major grid style={line width=.1pt,draw=gray!20},
        yminorticks=false,
        xticklabels={,,},
        title = {(a) Performance ratio, full accuracy}
        ]
        
        \addplot[clr7] table[x=tau-r , y=qics-rf] {\data}; 
        \addplot[clr1] table[x=tau-r , y=qics*-rf] {\data};
        \addplot[clr4] table[x=tau-r , y=cvxquad-rf] {\data};

    \end{semilogxaxis}

    \begin{semilogxaxis}[
        name=secondplot,
        at={(mainplot.north east)},
        xshift=1cm,
        anchor=north west,   
        width = 0.5\textwidth,
        height = 0.35\textwidth,
        xmin = 1e-2, xmax = 3600,
        ymin = 0, ymax=1,
        label style = {font=\footnotesize},
        ticklabel style = {font=\footnotesize},
        title style = {font=\footnotesize},
        no markers,
        grid=major,
        major grid style={line width=.1pt,draw=gray!20},
        yminorticks=false,
        yticklabels={,,},
        xticklabels={,,},
        title = {(b) Solution time profile, full accuracy}
        ]
        
        \addplot[clr7] table[x=tau-t , y=qics-tf] {\data}; 
        \addplot[clr1] table[x=tau-t , y=qics*-tf] {\data};
        \addplot[clr4] table[x=tau-t , y=cvxquad-tf] {\data};
        
    \end{semilogxaxis}

    \begin{semilogxaxis}[
        name=thirdplot,
        at={(mainplot.south west)},
        yshift=-1cm,
        anchor=north west,   
        width = 0.5\textwidth,
        height = 0.35\textwidth,
        xmin = 1, xmax = 10000,
        ymin = 0, ymax=1,
        xlabel=Relative performance ratio,
        ylabel=Ratio of problems solved,
        label style = {font=\footnotesize},
        ticklabel style = {font=\footnotesize},
        title style = {font=\footnotesize},
        no markers,
        grid=major,
        major grid style={line width=.1pt,draw=gray!20},
        yminorticks=false,
        title = {(c) Performance ratio, low accuracy}
        ]
        
        \addplot[clr7] table[x=tau-r , y=qics-rn] {\data}; 
        \addplot[clr1] table[x=tau-r , y=qics*-rn] {\data};
        \addplot[clr4] table[x=tau-r , y=cvxquad-rn] {\data};
        
    \end{semilogxaxis}

    \begin{semilogxaxis}[
        at={(thirdplot.north east)},
        xshift=1cm,
        anchor=north west,   
        width = 0.5\textwidth,
        height = 0.35\textwidth,
        xmin = 1e-2, xmax = 3600,
        ymin = 0, ymax=1,
        xlabel=Solve time (s),
        label style = {font=\footnotesize},
        ticklabel style = {font=\footnotesize},
        title style = {font=\footnotesize},
        no markers,
        grid=major,
        major grid style={line width=.1pt,draw=gray!20},
        yminorticks=false,
        yticklabels={,,},
        title = {(d) Solution time profile, low accuracy},
        legend cell align={left},
        legend columns=-1,
        legend style={/tikz/every even column/.append style={column sep=0.25cm}, font=\footnotesize, at={(0.35,-0.3)},cells={line width=1.5pt}},
        ]
        
        \addplot[clr7] table[x=tau-t , y=qics-tn] {\data}; 
        \addplot[clr1] table[x=tau-t , y=qics*-tn] {\data};
        \addplot[clr4] table[x=tau-t , y=cvxquad-tn] {\data};

        \legend{
            QICS,
            QICS*,
            CVXQUAD
        }
        
    \end{semilogxaxis}

\end{tikzpicture}

%% file: main.bbl
\begin{thebibliography}{10}
\providecommand{\url}[1]{#1}
\csname url@samestyle\endcsname
\providecommand{\newblock}{\relax}
\providecommand{\bibinfo}[2]{#2}
\providecommand{\BIBentrySTDinterwordspacing}{\spaceskip=0pt\relax}
\providecommand{\BIBentryALTinterwordstretchfactor}{4}
\providecommand{\BIBentryALTinterwordspacing}{\spaceskip=\fontdimen2\font plus
\BIBentryALTinterwordstretchfactor\fontdimen3\font minus \fontdimen4\font\relax}
\providecommand{\BIBforeignlanguage}[2]{{%
\expandafter\ifx\csname l@#1\endcsname\relax
\typeout{** WARNING: IEEEtran.bst: No hyphenation pattern has been}%
\typeout{** loaded for the language `#1'. Using the pattern for}%
\typeout{** the default language instead.}%
\else
\language=\csname l@#1\endcsname
\fi
#2}}
\providecommand{\BIBdecl}{\relax}
\BIBdecl

\bibitem{fawzi2023optimal}
H.~Fawzi and J.~Saunderson, ``Optimal self-concordant barriers for quantum relative entropies,'' \emph{SIAM Journal on Optimization}, vol.~33, no.~4, pp. 2858--2884, 2023.

\bibitem{fujii1989relative}
J.~I. Fujii and E.~Kamei, ``Relative operator entropy in noncommutative information theory,'' \emph{Mathematica Japonica}, vol.~34, pp. 341--348, 1989.

\bibitem{kubo1980means}
F.~Kubo and T.~Ando, ``Means of positive linear operators,'' \emph{Mathematische Annalen}, vol. 246, pp. 205--224, 1980.

\bibitem{nesterov2012towards}
Y.~Nesterov, ``Towards non-symmetric conic optimization,'' \emph{Optimization methods and software}, vol.~27, no. 4-5, pp. 893--917, 2012.

\bibitem{coey2023performance}
C.~Coey, L.~Kapelevich, and J.~P. Vielma, ``Performance enhancements for a generic conic interior point algorithm,'' \emph{Mathematical Programming Computation}, vol.~15, no.~1, pp. 53--101, 2023.

\bibitem{tunccel2001generalization}
L.~Tun{\c{c}}el, ``Generalization of primal—dual interior-point methods to convex optimization problems in conic form,'' \emph{Foundations of computational mathematics}, vol.~1, pp. 229--254, 2001.

\bibitem{myklebust2014interior}
T.~Myklebust and L.~Tun{\c{c}}el, ``Interior-point algorithms for convex optimization based on primal-dual metrics,'' \emph{arXiv preprint arXiv:1411.2129}, 2014.

\bibitem{myklebust2015primal}
T.~G. J.~J. Myklebust, ``On primal-dual interior-point algorithms for convex optimisation,'' Ph.D. dissertation, University of Waterloo, 2015.

\bibitem{nemirovski2005cone}
A.~Nemirovski and L.~Tun{\c{c}}el, ``“{C}one-free” primal-dual path-following and potential-reduction polynomial time interior-point methods,'' \emph{Mathematical Programming}, vol. 102, pp. 261--294, 2005.

\bibitem{karimi2020primal}
M.~Karimi and L.~Tun{\c{c}}el, ``Primal--dual interior-point methods for domain-driven formulations,'' \emph{Mathematics of Operations Research}, vol.~45, no.~2, pp. 591--621, 2020.

\bibitem{he2024exploiting}
K.~He, J.~Saunderson, and H.~Fawzi, ``Exploiting structure in quantum relative entropy programs,'' \emph{arXiv preprint arXiv:2407.00241}, 2024.

\bibitem{karimi2023efficient}
M.~Karimi and L.~Tuncel, ``Efficient implementation of interior-point methods for quantum relative entropy,'' \emph{arXiv preprint arXiv:2312.07438}, 2023.

\bibitem{karimi2024domain}
M.~Karimi and L.~Tun{\c{c}}el, ``Domain-driven solver ({DDS}) version 2.1: a {MATLAB}-based software package for convex optimization problems in domain-driven form,'' \emph{Mathematical Programming Computation}, vol.~16, no.~1, pp. 37--92, 2024.

\bibitem{fawzi2019semidefinite}
H.~Fawzi, J.~Saunderson, and P.~A. Parrilo, ``Semidefinite approximations of the matrix logarithm,'' \emph{Foundations of Computational Mathematics}, vol.~19, pp. 259--296, 2019.

\bibitem{mosek}
\BIBentryALTinterwordspacing
M.~ApS, \emph{The MOSEK optimization toolbox for MATLAB manual. Version 9.0.}, 2019. [Online]. Available: \url{http://docs.mosek.com/9.0/toolbox/index.html}
\BIBentrySTDinterwordspacing

\bibitem{vandenberghe2010cvxopt}
L.~Vandenberghe, ``The {CVXOPT} linear and quadratic cone program solvers,'' \emph{Online: http://cvxopt.org/documentation/coneprog.pdf}, 2010.

\bibitem{sturm1999using}
J.~F. Sturm, ``Using {SeDuMi} 1.02, a {MATLAB} toolbox for optimization over symmetric cones,'' \emph{Optimization methods and software}, vol.~11, no. 1-4, pp. 625--653, 1999.

\bibitem{nesterov1997self}
Y.~E. Nesterov and M.~J. Todd, ``Self-scaled barriers and interior-point methods for convex programming,'' \emph{Mathematics of Operations research}, vol.~22, no.~1, pp. 1--42, 1997.

\bibitem{nesterov1998primal}
------, ``Primal-dual interior-point methods for self-scaled cones,'' \emph{SIAM Journal on optimization}, vol.~8, no.~2, pp. 324--364, 1998.

\bibitem{fujisawa1997exploiting}
K.~Fujisawa, M.~Kojima, and K.~Nakata, ``Exploiting sparsity in primal-dual interior-point methods for semidefinite programming,'' \emph{Mathematical Programming}, vol.~79, pp. 235--253, 1997.

\bibitem{yamashita2010high}
M.~Yamashita, K.~Fujisawa, K.~Nakata, M.~Nakata, M.~Fukuda, K.~Kobayashi, and K.~Goto, ``A high-performance software package for semidefinite programs: {SDPA} 7,'' \emph{Tokyo, Japan}, 2010.

\bibitem{yamashita2012latest}
M.~Yamashita, K.~Fujisawa, M.~Fukuda, K.~Kobayashi, K.~Nakata, and M.~Nakata, ``Latest developments in the {SDPA} family for solving large-scale {SDPs},'' \emph{Handbook on semidefinite, conic and polynomial optimization}, pp. 687--713, 2012.

\bibitem{toh1999sdpt3}
K.-C. Toh, M.~J. Todd, and R.~H. T{\"u}t{\"u}nc{\"u}, ``{SDPT3}—a {MATLAB} software package for semidefinite programming, version 1.3,'' \emph{Optimization methods and software}, vol.~11, no. 1-4, pp. 545--581, 1999.

\bibitem{goemans2001approximation}
M.~X. Goemans and D.~Williamson, ``Approximation algorithms for {MAX-3-CUT} and other problems via complex semidefinite programming,'' in \emph{Proceedings of the thirty-third annual ACM symposium on Theory of computing}, 2001, pp. 443--452.

\bibitem{fujisawa2002sdpa}
K.~Fujisawa, M.~Kojima, K.~Nakata, and M.~Yamashita, ``{SDPA} (semidefinite programming algorithm) user’s manual --- version 6.2.0,'' \emph{Department of Mathematical and Com-puting Sciences, Tokyo Institute of Technology. Research Reports on Mathematical and Computing Sciences Series B: Operations Research}, vol.~5, p.~6, 2002.

\bibitem{friberg2016cblib}
H.~A. Friberg, ``{CBLIB} 2014: a benchmark library for conic mixed-integer and continuous optimization,'' \emph{Mathematical Programming Computation}, vol.~8, pp. 191--214, 2016.

\bibitem{picos}
G.~Sagnol and M.~Stahlberg, ``{PICOS}: A {Python} interface to conic optimization solvers,'' \emph{Journal of Open Source Software}, vol.~7, no.~70, p. 3915, Feb. 2022.

\bibitem{coles2016numerical}
P.~J. Coles, E.~M. Metodiev, and N.~L{\"u}tkenhaus, ``Numerical approach for unstructured quantum key distribution,'' \emph{Nature communications}, vol.~7, no.~1, p. 11712, 2016.

\bibitem{winick2018reliable}
A.~Winick, N.~L{\"u}tkenhaus, and P.~J. Coles, ``Reliable numerical key rates for quantum key distribution,'' \emph{Quantum}, vol.~2, p.~77, 2018.

\bibitem{hu2022robust}
H.~Hu, J.~Im, J.~Lin, N.~L{\"u}tkenhaus, and H.~Wolkowicz, ``Robust interior point method for quantum key distribution rate computation,'' \emph{Quantum}, vol.~6, p. 792, 2022.

\bibitem{nesterov2018lectures}
Y.~Nesterov, \emph{Lectures on Convex Optimization}.\hskip 1em plus 0.5em minus 0.4em\relax Springer, 2018.

\bibitem{lorente2024quantum}
A.~G. Lorente, P.~V. Parellada, M.~Castillo-Celeita, and M.~Ara{\'u}jo, ``Quantum key distribution rates from non-symmetric conic optimization,'' \emph{arXiv preprint arXiv:2407.00152}, 2024.

\bibitem{skajaa2015homogeneous}
A.~Skajaa and Y.~Ye, ``A homogeneous interior-point algorithm for nonsymmetric convex conic optimization,'' \emph{Mathematical Programming}, vol. 150, pp. 391--422, 2015.

\bibitem{nesterov1994interior}
Y.~Nesterov and A.~Nemirovskii, \emph{Interior-point polynomial algorithms in convex programming}.\hskip 1em plus 0.5em minus 0.4em\relax SIAM, 1994.

\bibitem{papp2017homogeneous}
D.~Papp and S.~Y{\i}ld{\i}z, ``On "a homogeneous interior-point algorithm for non-symmetric convex conic optimization",'' \emph{arXiv preprint arXiv:1712.00492}, 2017.

\bibitem{hiai2014introduction}
F.~Hiai and D.~Petz, \emph{Introduction to matrix analysis and applications}.\hskip 1em plus 0.5em minus 0.4em\relax Springer Science \& Business Media, 2014.

\bibitem{higham2008functions}
N.~J. Higham, \emph{Functions of matrices: theory and computation}.\hskip 1em plus 0.5em minus 0.4em\relax SIAM, 2008.

\bibitem{coey2023conic}
C.~Coey, L.~Kapelevich, and J.~P. Vielma, ``Conic optimization with spectral functions on {E}uclidean {J}ordan algebras,'' \emph{Mathematics of Operations Research}, vol.~48, no.~4, pp. 1906--1933, 2023.

\bibitem{faybusovich2020self}
L.~Faybusovich and C.~Zhou, ``Self-concordance and matrix monotonicity with applications to quantum entanglement problems,'' \emph{Applied Mathematics and Computation}, vol. 375, p. 125071, 2020.

\bibitem{hiai2017different}
F.~Hiai and M.~Mosonyi, ``Different quantum f-divergences and the reversibility of quantum operations,'' \emph{Reviews in Mathematical Physics}, vol.~29, no.~07, p. 1750023, 2017.

\bibitem{bhatia2013matrix}
R.~Bhatia, \emph{Matrix analysis}.\hskip 1em plus 0.5em minus 0.4em\relax Springer Science \& Business Media, 2013, vol. 169.

\bibitem{mittelmann}
H.~Mittelmann, ``Benchmarks for optimization software,'' \url{https://plato.asu.edu/bench.html}.

\bibitem{goulart2024clarabel}
P.~J. Goulart and Y.~Chen, ``Clarabel: An interior-point solver for conic programs with quadratic objectives,'' \emph{arXiv preprint arXiv:2405.12762}, 2024.

\bibitem{dolan2002benchmarking}
E.~D. Dolan and J.~J. Mor{\'e}, ``Benchmarking optimization software with performance profiles,'' \emph{Mathematical programming}, vol.~91, pp. 201--213, 2002.

\bibitem{borchers1999sdplib}
B.~Borchers, ``{SDPLIB} 1.2, a library of semidefinite programming test problems,'' \emph{Optimization Methods and Software}, vol.~11, no. 1-4, pp. 683--690, 1999.

\bibitem{skrzypczyk2023semidefinite}
P.~Skrzypczyk and D.~Cavalcanti, \emph{Semidefinite Programming in Quantum Information Science}.\hskip 1em plus 0.5em minus 0.4em\relax IOP Publishing, 2023.

\bibitem{watrous2012simpler}
J.~Watrous, ``Simpler semidefinite programs for completely bounded norms,'' \emph{arXiv preprint arXiv:1207.5726}, 2012.

\bibitem{cole2023quantum}
S.~Cole, M.~Eckstein, S.~Friedland, and K.~{\.Z}yczkowski, ``On quantum optimal transport,'' \emph{Mathematical Physics, Analysis and Geometry}, vol.~26, no.~2, p.~14, 2023.

\bibitem{doherty2004complete}
A.~C. Doherty, P.~A. Parrilo, and F.~M. Spedalieri, ``Complete family of separability criteria,'' \emph{Physical Review A}, vol.~69, no.~2, p. 022308, 2004.

\bibitem{siddhu2022five}
V.~Siddhu and S.~Tayur, ``Five starter pieces: Quantum information science via semidefinite programs,'' in \emph{Tutorials in Operations Research: Emerging and Impactful Topics in Operations}.\hskip 1em plus 0.5em minus 0.4em\relax INFORMS, 2022, pp. 59--92.

\bibitem{navascues2007bounding}
M.~Navascu{\'e}s, S.~Pironio, and A.~Ac{\'\i}n, ``Bounding the set of quantum correlations,'' \emph{Physical Review Letters}, vol.~98, no.~1, p. 010401, 2007.

\bibitem{navascues2008convergent}
------, ``A convergent hierarchy of semidefinite programs characterizing the set of quantum correlations,'' \emph{New Journal of Physics}, vol.~10, no.~7, p. 073013, 2008.

\bibitem{garner2024introducing}
A.~J. Garner and M.~Ara{\'u}jo, ``Introducing {Moment}: A toolkit for semi-definite programming with moment matrices,'' \emph{arXiv preprint arXiv:2406.15559}, 2024.

\bibitem{brown2021device}
P.~Brown, H.~Fawzi, and O.~Fawzi, ``Device-independent lower bounds on the conditional von {N}eumann entropy,'' \emph{Quantum}, vol.~8, p. 1445, 2024.

\bibitem{nakata2001variational}
M.~Nakata, H.~Nakatsuji, M.~Ehara, M.~Fukuda, K.~Nakata, and K.~Fujisawa, ``Variational calculations of fermion second-order reduced density matrices by semidefinite programming algorithm,'' \emph{The Journal of Chemical Physics}, vol. 114, no.~19, pp. 8282--8292, 2001.

\bibitem{nakata2002density}
M.~Nakata, M.~Ehara, and H.~Nakatsuji, ``Density matrix variational theory: Application to the potential energy surfaces and strongly correlated systems,'' \emph{The Journal of Chemical Physics}, vol. 116, no.~13, pp. 5432--5439, 2002.

\bibitem{yamashita2003implementation}
M.~Yamashita, K.~Fujisawa, and M.~Kojima, ``Implementation and evaluation of {SDPA} 6.0 (semidefinite programming algorithm 6.0),'' \emph{Optimization Methods and Software}, vol.~18, no.~4, pp. 491--505, 2003.

\bibitem{aps2020mosek}
M.~ApS, ``{MOSEK} modeling cookbook,'' 2020.

\bibitem{horodecki1996separability}
M.~Horodecki, P.~Horodecki, and R.~Horodecki, ``Separability of mixed states: necessary and sufficient conditions,'' \emph{Physics Letters A}, vol. 223, no.~1, pp. 1--8, 1996.

\bibitem{fawzi2018efficient}
H.~Fawzi and O.~Fawzi, ``Efficient optimization of the quantum relative entropy,'' \emph{Journal of Physics A: Mathematical and Theoretical}, vol.~51, no.~15, p. 154003, 2018.

\bibitem{he2024efficient}
K.~He, J.~Saunderson, and H.~Fawzi, ``Efficient computation of the quantum rate-distortion function,'' \emph{Quantum}, vol.~8, p. 1314, Apr. 2024.

\bibitem{fawzi2023entropy}
H.~Fawzi, O.~Fawzi, and S.~O. Scalet, ``Entropy constraints for ground energy optimization,'' \emph{Journal of Mathematical Physics}, vol.~65, no.~3, 2024.

\bibitem{bertsimas2023new}
D.~Bertsimas, R.~Cory-Wright, and J.~Pauphilet, ``A new perspective on low-rank optimization,'' \emph{Mathematical Programming}, vol. 202, no.~1, pp. 47--92, 2023.

\bibitem{berta2017variational}
M.~Berta, O.~Fawzi, and M.~Tomamichel, ``On variational expressions for quantum relative entropies,'' \emph{Letters in Mathematical Physics}, vol. 107, pp. 2239--2265, 2017.

\bibitem{huang2024semi}
Z.~Huang and M.~M. Wilde, ``Semi-definite optimization of the measured relative entropies of quantum states and channels,'' \emph{arXiv preprint arXiv:2406.19060}, 2024.

\bibitem{fawzi2024certified}
H.~Fawzi, O.~Fawzi, and S.~O. Scalet, ``Certified algorithms for equilibrium states of local quantum {H}amiltonians,'' \emph{Nature Communications}, vol.~15, no.~1, p. 7394, 2024.

\bibitem{andersen1995presolving}
E.~D. Andersen and K.~D. Andersen, ``Presolving in linear programming,'' \emph{Mathematical programming}, vol.~71, pp. 221--245, 1995.

\bibitem{drusvyatskiy2017many}
D.~Drusvyatskiy, H.~Wolkowicz \emph{et~al.}, ``The many faces of degeneracy in conic optimization,'' \emph{Foundations and Trends{\textregistered} in Optimization}, vol.~3, no.~2, pp. 77--170, 2017.

\bibitem{he2024bregman}
K.~He, J.~Saunderson, and H.~Fawzi, ``A {B}regman proximal perspective on classical and quantum {B}lahut-{A}rimoto algorithms,'' \emph{IEEE Transactions on Information Theory}, vol.~70, no.~8, pp. 5710--5730, 2024.

\bibitem{wolkowicz2004solving}
H.~Wolkowicz, ``Solving semidefinite programs using preconditioned conjugate gradients,'' \emph{Optimization Methods and Software}, vol.~19, no.~6, pp. 653--672, 2004.

\bibitem{tang2024feasible}
T.~Tang and K.-C. Toh, ``A feasible method for general convex low-rank {SDP} problems,'' \emph{SIAM Journal on Optimization}, vol.~34, no.~3, pp. 2169--2200, 2024.

\end{thebibliography}
